\documentclass[gdplain]{geradwp}

\PassOptionsToPackage{monochrome}{xcolor}
\PassOptionsToPackage{hyphens}{url}

\usepackage{lipsum}

\usepackage[english]{babel}
\usepackage{lastpage}
\usepackage{xspace}
\usepackage{fancyhdr}
\usepackage{geometry}
\usepackage{setspace}
\usepackage{array}
\usepackage{multirow}
\usepackage[utf8]{inputenc}
\usepackage[T1]{fontenc}
\usepackage{csquotes}
\usepackage{lmodern}
\usepackage[dvipsnames]{xcolor}
\usepackage{hyperref}
\usepackage{caption}
\usepackage{enumitem}
\usepackage{graphicx}
\usepackage{float}
\usepackage{amsmath}
\usepackage{amsthm}
\usepackage{dsfont}
\usepackage{stmaryrd}
\usepackage[b]{esvect}
\usepackage[ruled]{algorithm}
\usepackage[noend]{algpseudocode}
\usepackage[style=numeric,maxbibnames=10,backend=bibtex]{biblatex}
\GDcoverpagewhitespace{6.8cm}
\hypersetup{colorlinks,%
citecolor={blue}, 
urlcolor={blue},
linkcolor={blue},
breaklinks={true}
}

\makeatletter
\ifthenelse{\isundefined{\ALG@name}}{}%
{%
\renewcommand{\ALG@name}{\sffamily\footnotesize Algorithm}
}
\makeatother
\ifthenelse{\isundefined{\SetAlCapNameFnt}}{}%
{%
\SetAlCapNameFnt{\footnotesize}
\SetAlCapFnt{\sffamily\footnotesize}
}

\newcommand{\guillemetseng}[1]{\textquotedblleft#1\textquotedblright}
\newcommand{\inlinequote}[1]{\guillemetseng{#1}}

\newcolumntype{L}[1]{>{\raggedright\arraybackslash}p{#1}}
\newcolumntype{R}[1]{>{\raggedleft\arraybackslash}p{#1}}
\newcolumntype{M}[1]{>{\centering\arraybackslash}m{#1}}

\newlength\aetmplength

\newcommand{\Acal}{\mathcal{A}}\newcommand{\Abf}{\mathbf{A}}
\newcommand{\Bcal}{\mathcal{B}}
\newcommand{\Ccal}{\mathcal{C}}\newcommand{\Crm}{\mathrm{C}}\newcommand{\Cbb}{\mathbb{C}}
\newcommand{\Dcal}{\mathcal{D}}

\newcommand{\Hcal}{\mathcal{H}}
\newcommand{\Ibb}{\mathbb{I}}
\newcommand{\Jbb}{\mathbb{J}}

\newcommand{\Lcal}{\mathcal{L}}
\newcommand{\Mbb}{\mathbb{M}}
\newcommand{\Nbb}{\mathbb{N}}
\newcommand{\Obb}{\mathbb{O}}
\newcommand{\Pbf}{\mathbf{P}}\newcommand{\Pbb}{\mathbb{P}}

\newcommand{\Rbb}{\mathbb{R}}
\newcommand{\Srm}{\mathrm{S}}
\newcommand{\Tcal}{\mathcal{T}}

\newcommand{\Wcal}{\mathcal{W}}\newcommand{\Wbb}{\mathbb{W}}
\newcommand{\Xbb}{\mathbb{X}}
\newcommand{\Ybb}{\mathbb{Y}}

\newcommand{\algoname}[1]{\texttt{#1}\xspace}
\newcommand{\dfo}{\algoname{DFO}}

\newcommand{\dsm}{\algoname{dsm}}
\newcommand{\cdsm}{\algoname{cdsm}}

\newcommand{\bfgs}{\algoname{bfgs}}

\newcommand{\searchstep}{\algoname{search}}
\newcommand{\pollstep}{\algoname{poll}}
\newcommand{\coveringstep}{\algoname{covering}}

\newcommand{\failurestep}{\algoname{failure}}
\newcommand{\attackstep}{\algoname{attack}}
\newcommand{\attackset}{\Acal_\Lcal}
\newcommand{\attacksetfeas}{\attackset^\mathrm{f}}
\newcommand{\attacksetsuccess}{\attackset^+}
\newcommand{\attacksetopt}{\attackset^*}

\newcommand{\nn}{\algoname{NN}}
\newcommand{\nns}{{\nn}s\xspace}
\newcommand{\pinn}{\algoname{PINN}}
\newcommand{\resnet}{\algoname{ResNet18}}
\newcommand{\fgsm}{\algoname{fgsm}}
\newcommand{\ffgsm}{\algoname{ffgsm}}
\newcommand{\rfgsm}{\algoname{rfgsm}}
\newcommand{\pgd}{\algoname{pgd}}
\newcommand{\relu}{\algoname{ReLU}}
\newcommand{\atck}{\algoname{atck}}
\newcommand{\hybr}{\algoname{hybr}}
\newcommand{\bprp}{\algoname{back}}

 \newcommand{\TcalS}{\Tcal_{\tt{S}}} \newcommand{\tS}{t_{\tt{S}}} 
 \newcommand{\TcalC}{\Tcal_{\tt{C}}} \newcommand{\tC}{t_{\tt{C}}} 
 \newcommand{\TcalP}{\Tcal_{\tt{P}}} \newcommand{\tP}{t_{\tt{P}}} 
\newcommand{\DcalA}{\Dcal_{\tt{A}}} \newcommand{\TcalA}{\Tcal_{\tt{A}}} \newcommand{\tA}{t_{\tt{A}}}

\newcommand{\abcrown}{$(\alpha,\beta)$-\texttt{CROWN}\xspace}
\newcommand{\squarerrorloss}{\Lcal_{\mathrm{SE}}}
\newcommand{\crossentropyloss}{\Lcal_{\mathrm{CE}}}
\newcommand{\wslp}{\texttt{WSLP}\xspace}
\newcommand{\bary}{\algoname{bary}}
\newcommand{\softmax}{\algoname{softmax}}
\newcommand{\backprop}{\algoname{backprop}}
\newcommand{\warcraft}{\algoname{warcraft}}
\newcommand{\costmap}{\algoname{costmap}}
\newcommand{\pathmap}{\mathrm{p}}
\newcommand{\costpath}{\algoname{costpath}}
\newcommand{\counterfactual}{\mathrm{cfa}}

\newcommand{\Wcalini}{\Wcal_\mathrm{ini}}
\newcommand{\pathmapini}{\pathmap^*_\mathrm{ini}}
\newcommand{\xini}{x_\mathrm{ini}}
\newcommand{\Ccalini}{\Ccal_\mathrm{ini}}

\renewcommand{\algorithmicfor}{\textbf{for}\xspace}

\renewcommand{\algorithmicfunction}{\textbf{function}\xspace}
\renewcommand{\algorithmicreturn}{\textbf{return}\xspace}
\algrenewcommand{\textproc}{\emph}

\newcommand{\defequal}{\triangleq}
\newcommand{\dsum}[2]{\sum\limits_{#1}^{#2}}

\DeclareMathOperator*{\argmin}{argmin}
\DeclareMathOperator*{\argmax}{argmax}
\DeclareMathOperator*{\minimize}{minimize}
\DeclareMathOperator*{\maximize}{maximize}
\newcommand{\1}{\mathds{1}}

\newcommand{\norm}[1]{\left\lVert{#1}\right\rVert}
\newcommand{\dotprod}[2]{\left<#1~,~#2\right>}
\newcommand{\textabs }[1]{\lvert{#1}\rvert}
\newcommand{\textnorm}[1]{\lVert{#1}\rVert}
\newcommand{\textdotprod}[2]{\left<#1,#2\right>}
\newcommand{\transpose}[1]{#1^\top}
\def\llb{\llbracket}
\def\rrb{\rrbracket}
\newcommand{\interior}{\mathrm{int}}
\newcommand{\cl}{\mathrm{cl}}
\newcommand{\dist}{\mathrm{dist}}

\newtheorem{assumption}{Assumption}
\newtheorem{proposition}{Proposition}

\newtheorem{definition}{Definition}
\newtheorem{theorem}{Theorem}
\newtheorem{remark}{Remark}

\usepackage{forloop}
\newcounter{ct}
\newcommand{\markdent}[1]{\hspace{-\algorithmicindent}\forloop{ct}{0}{\value{ct} < #1}{\hspace{1em}}}
\newcommand{\IState}[1]{\Statex\markdent{#1}}

\newcommand{\fct}[5]{
    \setlength{\arraycolsep}{2pt}
    #1:\left\{\begin{array}{ccl}
        #3 & \to     & #5\\
        #2 & \mapsto & #4
    \end{array}\right.}

\newcommand{\compactarray}[1]{\begingroup\scriptstyle\renewcommand{\arraystretch}{1.0}\begin{array}{c}#1\end{array}\endgroup}
\newcommand{\problemoptimfree}[3]{\underset{\compactarray{#2}}{#1} \quad #3}
\newcommand{\problemoptimoneline}[4]{\underset{\compactarray{#2}}{#1} \quad #3 \quad \mbox{subject to} \quad #4}
\newcommand{\problemoptim}[4]{
    \begin{array}{cl}
        \underset{\compactarray{#2}}{#1} & \begin{array}{l}#3\end{array} \\ \\
        \mbox{subject to} & \begin{array}[t]{ll}#4\end{array}
    \end{array}
}

\setlist[enumerate]{noitemsep, topsep=-\parskip}
\setlist[itemize]{noitemsep, topsep=-\parskip}

\GDtitle{Optimization by Directional Attacks: Solving Problems with Neural Network Surrogates}
\GDmonth{TBD}{TBD}
\GDyear{2025}
\GDnumber{TBD}
\GDauthorsShort{P.-Y. Bouchet, T. Vidal}
\GDauthorsCopyright{Bouchet, Vidal}
\GDpostpubcitation{Hamel, Benoit, Karine H\'ebert (2021). \og Un exemple de citation \fg, \emph{Journal of Journals}, vol. X issue Y, p. n-m}{https://www.gerad.ca/fr}
\GDrevised{Août}{August}{2025}

\begin{document}


\begin{GDtitlepage}

\begin{GDauthlist}
\GDauthitem{Pierre-Yves Bouchet \ref{affil:polymtl}\GDrefsep\ref{affil:cirrelt}}
\GDauthitem{Thibaut Vidal \ref{affil:polymtl}\GDrefsep\ref{affil:cirrelt}}
\end{GDauthlist}

\begin{GDaffillist}
\GDaffilitem{affil:polymtl}{École Polytechnique de Montréal, Montr\'eal (Qc), Canada, H3T 1J4}
\GDaffilitem{affil:cirrelt}{CIRRELT, Montr\'eal (Qc), Canada, H3T 1J4}
\end{GDaffillist}

\begin{GDemaillist}
\GDemailitem{pierre-yves.bouchet@polymtl.ca}
\GDemailitem{thibaut.vidal@polymtl.ca}
\end{GDemaillist}

\end{GDtitlepage}


\GDabstracts

\begin{GDabstract}{Abstract}
    This paper tackles optimization problems whose objective and constraints involve a trained \emph{Neural Network} (\nn), where the goal is to maximize~$f(\Phi(x))$ subject to~$c(\Phi(x)) \leq 0$, with~$f$ smooth,~$c$ general and non-stringent, and~$\Phi$ a \nn already trained and not given as a glass-box.
    We address two challenges regarding this problem: identifying ascent directions for local search, and ensuring reliable convergence towards relevant local solutions.
    To this end, we re-purpose the notion of \emph{directional \nn attacks} as efficient optimization subroutines, since directional \nn attacks are designed to efficiently compute perturbations of~$x$ that steer~$\Phi(x)$ in prescribed directions.
    Precisely, we develop an \attackstep operator that computes attacks of~$\Phi$ at any~$x$ along the direction~$\nabla f(\Phi(x))$.
    Then, we propose a hybrid algorithm combining the \attackstep operator with \emph{derivative-free optimization} (\dfo) techniques, designed for numerical reliability by remaining oblivious to the structure of the problem.
    We consider the \emph{Covering Direct Search Method (\cdsm)}, which offers asymptotic guarantees to converge to a local solution under mild assumptions on the problem.
    The resulting method alternates between attack-based steps for heuristic yet fast local improvements of the current incumbent solution and \cdsm steps for certified convergence and numerical reliability.
    Experiments on three problems show that this hybrid approach outperforms standard baselines.

\end{GDabstract}

\begin{GDacknowledgements}
We are grateful to Paul A. Patience and Bruno Blais, from the Department of Chemical Engineering at Polytechnique Montréal, for providing the Physics-Informed Neural Network used in our numerical experiment in Section~\ref{section:numerical/bio_pinn}.
This research was supported by SCALE-AI through its Research Chairs program.
\end{GDacknowledgements}


\clearpage\newpage
\GDarticlestart

\section{Introduction}
\label{section:introduction}

\emph{Neural Networks} (\nns) are modelling tools acknowledged across various domains.
Beyond classification tasks (such as in computer vision~\cite{ZeFaWeShJu22CNNsurvey} and medical diagnostics~\cite{SaKu22CNNMedicalSurvey}), \nns are increasingly deployed for applications such as physics-informed learning~\cite{CuDCGiRoRaPi22PINNsurvey,ShFeChZhCaJi25PINNSurvey} and uncertainty quantification~\cite{KaKhHoNa18NNuncertaintySurvey}.
This allows \nns to be used as \emph{digital twins}~\cite{TaXiQiChJi22DigitalTwin} of a physical system, that is, rich models of the system allowing for real-time decision-making.
This evolving use has given rise to a class of optimization problems in which a trained \nn~$\Phi: \Rbb^n \to \Rbb^m$ acts as a surrogate function mapping the decision variables to a system response.
The objective is to optimize the output of a task-specific goal function~$f: \Rbb^m \to \Rbb$ defined over the space of outputs of~$\Phi$, subject to constraints defined via a constraint function~$c: \Rbb^m \to \Rbb^p$.
This leads to what we refer to as the composite \emph{problem of optimization through a neural network},
\begin{equation}
    \label{problem:P}
    \tag{$\Pbf$}
    \problemoptimoneline
        {\maximize}
        {x \in \Rbb^n}
        {f(\Phi(x))}
        {c(\Phi(x)) \leq 0},
\end{equation}
where~$(n, m, p) \in (\Nbb^*)^3$ denote the dimensions of, respectively, the variables space, the \nn output space, and the constraints function output space.
We study Problem~\eqref{problem:P} under the following Assumption~\ref{assumption:problem}.

\begin{assumption}[Problem requirements]
    \label{assumption:problem}
    The \nn~$\Phi$ encodes a continuous function, the goal function~$f$ is differentiable with gradient~$\nabla f$, and the feasible set~$F \defequal \{x \in \Rbb^n: c(\Phi(x)) \leq 0\}$ induced by the constraints function~$c$ is nonempty, ample (that is,~$F \subseteq \cl(\interior(F))$) and compact.
\end{assumption}

We consider the \nn~$\Phi$ as already trained and not tuneable anymore.
Moreover, we consider that the mathematical function encoded by~$\Phi$ is not explicitly exposed.
We then work in a \emph{derivative-free optimization} (\dfo)~\cite{AuHa2017,CoScVibook,LaMeWi2019} setting.
In particular, we do not assume access to the architecture and weights of~$\Phi$.
Instead, we typically consider that~$\Phi$ is given in a compiled form that maximizes speed, at the cost of transparency.
Our main usecase is when~$\Phi$ retains \emph{backpropagations}~\cite[Section~6.5]{GoBeCo16DeepLearning} capability.
This requirement is met in most frameworks for \nns, as they may compile models allowing off-the-shelf backpropagations that do not expose the structure of the \nn.
Yet, our methodology is also theoretically compatible with cases where~$\Phi$ is given as a \emph{black-box}.

To our best knowledge, only a few papers tackle Problem~\eqref{problem:P}, and most do so under more restrictive settings.
They typically assume that~$f$ is linear and~$F$ is polyhedral, and moreover that~$\Phi$ is provided as a glass-box \nn~\cite{PeTs22NNsampling,PhReTaToSe25OptimTrainedSparseNN,PlHaKlGaSuSa25ReLUNNInOptim,ToCaSe24RelaxingWalk,WaLoCaBe23OptimEnsembleTrainedNN} (that is, the architecture and weights of the \nn are all given).
Yet, the widespread adoption of \nns across industrial and engineering applications suggests that instances of Problem~\eqref{problem:P} involving less permissive \nn access (e.g., \nns in compiled form) are increasingly common.
Let us present two contexts where this broader setting naturally appears.

\paragraph{Simulation-based optimization.}
\emph{Simulation-based optimization}~\cite{AlAuGhKoLed2020,ATAuDiGhLDLeLGTr25solar,BB_list} is a common \dfo context that maximizes~$f(B(x))$ subject to~$c(B(x)) \leq 0$ over~$x \in \Rbb^n$, where~$f$ and~$c$ are defined as in Problem~\eqref{problem:P} and the mapping~$B: \Rbb^n \to \Rbb^m$ denotes an intractable process that, typically, runs a costly numerical simulation.
This setting, popular in engineering design~\cite{LaMe24structureAware,LaMeZh21ManifoldComposite}, has the same composite structure as Problem~\eqref{problem:P}, differing only in that the intermediate mapping is the simulator~$B$ rather than the \nn~$\Phi$.
Yet, an increasing trend replaces simulators~$B$ by trained \nns~$\Phi$ acting as surrogates~\cite[Chapter~15]{AuHa2017}, which yields instances of Problem~\eqref{problem:P} with key advantages:~$\Phi$ is typically runs orders of magnitude faster to run that~$B$; and even if~$\Phi$ is given as a compiled file rather than a glass-box, its neural structure can be exploited.
This motivates the development of dedicated optimization algorithms.

\paragraph{Counterfactual explanations of \nns.}
Counterfactual explanations~\cite{VeBoHoHiDiSh22CounterfactualReview,WaMiRu18CounterfactualBB} for a \nn consist in minimal changes to a given input that induce a desired change in the output.
When considering a classification task for a \nn classifier~$\Phi: \Rbb^n \to \llb1,m\rrb$, counterfactual explanations~$x_\counterfactual \in \Rbb^n$ for an input~$\xini \in \Rbb^n$ and a target class~$y^\sharp \neq \Phi(\xini)$ are solutions to~\inlinequote{$\minimize_{x \in \Rbb^n} \textnorm{x-\xini}^2$ subject to~$\Phi(x) = y^\sharp$}.
The notion also extends to contextual optimization~\cite{BoPrLePa25PrimalDualAlgoContextualStochasticCO,SaChDeFoFrVi25SurveyContextualCO}, where the \nn~$\Phi: \Rbb^n \to \Rbb^m$ maps a so-called \emph{context}~$x \in \Rbb^n$ to parameters for a downstream optimization task~\inlinequote{$\minimize_{y \in \Ccal} g(\Phi(x),y)$}, involving some function~$g$ usually convex and some set~$\Ccal$ that is typically combinatorial, for which an optimal decision rule~$y^*(x) \in \argmin_{y \in \Ccal} g(\Phi(x), y)$ may be selected.
In this setting, a counterfactual explanation of a context~$\xini \in \Rbb^n$ consists in a context~$x_\counterfactual \in \Rbb^n$ close to~$\xini$ and such that a given target decision~$y^\sharp \in \Ccal$, usually proposed by an expert or a benchmark policy, is preferable to~$y^*(\xini)$ for the problem induced by~$\Phi(x_\counterfactual)$~\cite{FoPaVi23DataDrivenOptim,VAFoPaVi24CFOPT}.
This results in the problem of \emph{$\varepsilon$-relative counterfactual}, for any~$\varepsilon \in \Rbb_+$, which aligns with Problem~\eqref{problem:P} since it is given by
\begin{equation*}
    \problemoptimoneline
        {\minimize}
        {x \in \Rbb^n}
        {\norm{x-\xini}^2}
        {(1+\varepsilon)g(\Phi(x), y^\sharp) \leq g(\Phi(x),y^*(\xini))}.
\end{equation*}

\paragraph{Contributions.}
This paper addresses the growing need to solve increasingly challenging instances of Problem~\eqref{problem:P} with the \nn~$\Phi$ not provided as a glass-box.
Our main insight is that the core challenge of identifying ascent directions for~$f \circ \Phi$ can be alleviated through \emph{directional \nn attacks}.
Indeed, directional \nn attacks are designed to find a perturbation of the input of a \nn (not given as a glass-box) that steers its outputs in a desired direction.
This paper solves Problem~\eqref{problem:P} by iteratively computing directional attacks of~$\Phi$ at the incumbent solution~$x^k$ in the direction given by~$\nabla f(\Phi(x^k))$.
Moreover, this approach may be hybridized with existing optimization algorithms, acting as an additional step at all iterations to search for ascent directions.
In particular, we hybridize this so-called \attackstep step with \dfo algorithms designed for global search capability and numerical reliability in problems with nonconvexity and nondifferentiability.
Then, in this paper,
\begin{enumerate}
    \vspace{0.25em}\item
        we formalize the concept of \emph{directional \nn attacks}, and we show that they can be repurposed to compute ascent directions for Problem~\eqref{problem:P}.
        We highlight that standard \nn attack algorithms can be used off-the-shelf for this purpose, and we describe how to embed them in a local optimization framework.
        While powerful when successful, this approach is inherently local and sometimes prone to failure for various technical reasons.
        This contribution forms the content of Section~\ref{section:optimization_attack};
    \vspace{0.25em}\item
        we hybridize the above attack-based approach with the \emph{covering direct search method} (\cdsm)~\cite{AuBoBo24Covering}, a \dfo method that is guaranteed to asymptotically converge to a local solution to Problem~\eqref{problem:P} under Assumption~\ref{assumption:problem}.
        Our resulting hybrid algorithm combines, at each iteration, two concepts with natural synergy: an \attackstep step leveraging the internal structure of~$\Phi$ for fast local improvement, and steps from \cdsm that allow for globalization strategies and asymptotic convergence towards a local solution.
        This contribution is detailed in Section~\ref{section:optimization_hybrid};
    \vspace{0.25em}\item
        we conduct numerical experiments on three problems to assess our hybrid algorithm.
        The first problem acts as a proof of concept involving a deep and complex \nn, while the other two are drawn from the simulation-based optimization and counterfactual explanation contexts highlighted in our motivation.
        Our experiments show that the \attackstep step is not a reliable standalone method, but it nevertheless contributes significantly to the performance of our hybrid method.
        Overall, our hybrid method appears faster than \dfo baselines thanks to the \attackstep step, and as reliable thanks to the steps from \cdsm.
        Details are given in Section~\ref{section:numerical}.
    \vspace{0.25em}
\end{enumerate}

We leave a literature review to Section~\ref{section:related_literature} and a discussion foreseeing future work to Section~\ref{section:conclusion}.
We conclude this section by introducing some notation used throughout the paper.

\paragraph{Notation.}
We denote by~$\textdotprod{\cdot}{\cdot}$ the dot product in~$\Rbb^m$, by~$\textnorm{\cdot}$ some norm in~$\Rbb^n$, by~$(e_i)_{i=1}^{n}$ the vectors of the canonical basis of~$\Rbb^n$, and by~$\1 \defequal \sum_{i=1}^{n}e_i$ the vector of all ones in~$\Rbb^n$.
For all~$r \in \Rbb_+$, all~$x \in \Rbb^n$, and all~$y \in \Rbb^m$, we respectively denote by~$\Bcal^n_r(x)$ and~$\Bcal^m_r(y)$ the closed ball of radius~$r$ centred at~$x$ and the closed ball of radius~$r$ in~$\Rbb^m$ centred at~$y$.
We omit the parentheses for the balls centred at~$0$, that is,~$\Bcal_r^n \defequal \Bcal^n_r(0)$ and~$\Bcal_r^m \defequal \Bcal^m_r(0)$.
For all \nn~$\Psi: \Rbb^a \to \Rbb^b$ (with~$(a,b) \in \Nbb^* \times \Nbb^*$), all~$x \in \Rbb^a$ and all~$d \in \Rbb^b$, we denote by~$\backprop(\Psi,x,d) \in \Rbb^a$ the output of a backpropagation of~$d$ through~$\Psi$ at~$x$.
We denote by~$\relu: v \in \Rbb^a \mapsto [\max(v_i,0)]_{i=1}^{a} \in \Rbb^a$ the \emph{rectified linear unit} function and by~$\softmax: v \in \Rbb^a \mapsto [\exp(v_i) / \sum_{j=1}^{a}\exp(v_j)]_{i=1}^{a} \in [0,1]^a$ the \emph{softmax} function, regardless of the dimension~$a \in \Nbb^*$ of the input vector~$v$.
For all couples~$(x_1, x_2) \in \Rbb^n \times \Rbb^n$, we say that~$x_2$ \emph{dominates}~$x_1$ (with respect to Problem~\eqref{problem:P}) if~$c(\Phi(x_2)) \leq 0$ and~$f(\Phi(x_2)) > f(\Phi(x_1))$.
Similarly, for all couples~$(x,d) \in \Rbb^n \times \Rbb^n$, we say that~$d$ is an \emph{ascent direction emanating from~$x$} (with respect to Problem~\eqref{problem:P}) if~$x+d$ dominates~$x$.
Note that we allow for non-unitary directions.

\section{Related Literature}
\label{section:related_literature}

Prior studies only tackle Problem~\eqref{problem:P} under more restrictive assumptions~\cite{PeTs22NNsampling,PhReTaToSe25OptimTrainedSparseNN,PlHaKlGaSuSa25ReLUNNInOptim,ToCaSe24RelaxingWalk,WaLoCaBe23OptimEnsembleTrainedNN}.
They typically assume that~$f$ is linear and~$F$ is polyhedral, and that~$\Phi$ is provided as a glass-box \nn.
These approaches effectively exploit the structure of~$\Phi$, but their requirements contrast with the reality of many practical applications, where~$\Phi$ may be deployed as an opaque or compiled artifact to maximize inference speed, which makes glass-box assumptions and layer-wise manipulations impractical.

Our setting relates to \dfo~\cite{AuHa2017,CoScVibook}, because of the potential nonsmoothness of~$f \circ \Phi$ and since~$\Phi$ is not a glass-box, which result in little exploitable structure on Problem~\eqref{problem:P}.
We are not aware of \dfo methods specifically designed for Problem~\eqref{problem:P}.
Then, we build our hybrid method upon the \emph{covering direct search method} (\cdsm)~\cite{AuBoBo24Covering} to inherit from its convergence towards a local solution to Problem~\eqref{problem:P}.
However, this choice is nonstringent and our core idea to hybridize attack techniques with a dedicated algorithm is compatible with many of the state-of-the-art \dfo algorithms surveyed in~\cite{LaMeWi2019}.

Our approach leveraging the neural structure of~$\Phi$ to identify ascent directions for~$f \circ \Phi$ relates to \nn attacks~\cite{SzZaSuBrErGoFe14NNAttackFirst}.
Algorithms computing \nn attacks cover a wide range of \nn accesses.
We refer to~\cite{zhang22babattack} for attacks of glass-box \nns, to \fgsm~\cite{Goodfellow15FGSM} and \pgd~\cite{Madry19PGD} and the many algorithms implemented in several solvers to compute \nn attacks~\cite{AdvertTorch,kim2020torchattacks,AdvRobustToolbox,CleverHans,rauber2017foolbox,zhang22babattack} for \nns providing accesses to backpropagation, and to algorithms such as~\cite{AlShChZhHsSr19GenAttack,AnCrFlHe20SquareAttack,GuGaYoWiWe19SimpleBlackBoxAttack} for \nns with black-box access.
Two surveys~\cite{SurveyTransferabilityAdversarialExamples,LiChHsLe21ReviewAdversarialAttackDefense} give a broader view of existing algorithms.
We design our \attackstep operator so that it may leverage any of them.
The notion of \nn attacks finally has connection with those of \emph{\nn verification}~\cite{KoBoHoRi24VerifNN,LiArLaStBaKo21VerifNN,Wu24VerifNNRobust}, another active research area~\cite{BrBaJoWu24VNN-COMP-2024,UbMi21ReviewVerification}.
However, the computational burden of a \nn verification is usually orders of magnitude greater than those of a \nn attack.
Most parts of the \abcrown solver~\cite{LiArLaStBaKo21VerifNN} address precisely this verification task, albeit primarily in binary classification contexts.

Finally, we stress that this work considers that the \nn~$\Phi$ is already trained.
Yet, the upstream training of~$\Phi$ may be optimized so that the resulting Problem~\eqref{problem:P} is easier to solve.
When~$\Phi$ models physical systems, conservation laws~\cite{HaMaAlGuMa23LearningPhysicalModels} or partial differential equations~\cite{NeMaKr23LearningDifferentiableSolvers} may be enforced by design, and constraints in a more generic sense may be satisfied in probability~\cite{UtMaMaWa25LearningHardConstraints}.
Similarly,~$\Phi$ may be trained so that the optimality conditions of the downstream problem are simpler to express~\cite{Ts21SobolevTrainedNNinOptim}.
When applicable, training~$\Phi$ accordingly simplifies the expression of Problem~\eqref{problem:P} or its solution process.

In short, our work significantly departs from the existing literature on optimization through \nns by addressing Problem~\eqref{problem:P} without assuming that~$\Phi$ is a glass-box \nn.
It also departs from the existing literature on \dfo by proposing a hybrid method that explicitly leverages the neural structure of~$\Phi$.

\section{Optimization Leveraging Directional \nn Attacks}
\label{section:optimization_attack}

This section addresses the first objective of the paper: establishing directional \nn attacks (formalized in Section~\ref{section:optimization_attack/notion}) as a practical tool for solving Problem~\eqref{problem:P}.
Specifically, in Section~\ref{section:optimization_attack/use_in_optim} we demonstrate that such attacks, when appropriately constructed, can be used to identify ascent directions.

\subsection{Directional \nn Attacks}
\label{section:optimization_attack/notion}

First, we contextualize the notion of \nn attack in Section~\ref{section:optimization_attack/notion/origin_illustration}.
Then we formalize directional \nn attacks in Section~\ref{section:optimization_attack/notion/definition}.
Finally, we discuss practical approaches for computing such attacks in Section~\ref{section:optimization_attack/notion/numerical_tools}.

\subsubsection{Background}
\label{section:optimization_attack/notion/origin_illustration}

The notion of \nn attacks originates from~\cite{SzZaSuBrErGoFe14NNAttackFirst}, which empirically demonstrates a striking vulnerability of many neural networks: small, carefully crafted perturbations of an input can lead to large changes in the output of the network.
The notion of \nn attacks was initially formalized for \nns focusing on classification as a central topic in adversarial machine learning.
A \emph{\nn attack} consists in finding an alteration of an input (within a ball of pre-defined radius) that changes the classification from the class predicted initially to any other.
Similarly, a \emph{targeted \nn attack} consists in altering the input in order to make it classified as a specific class.
In our context in this work, the \nns we consider may not perform classification, but nevertheless the notion of \nn attacks admits a direct extension that we formalize in Section~\ref{section:optimization_attack/notion/definition}.
Before going to this point, let us illustrate the observation from~\cite{SzZaSuBrErGoFe14NNAttackFirst} that a \nn classification may be drastically vulnerable to slight alterations of the input.

\begin{figure}[!ht]
    \centering
    \includegraphics[width=\linewidth, trim= 0 115 0 0, clip]{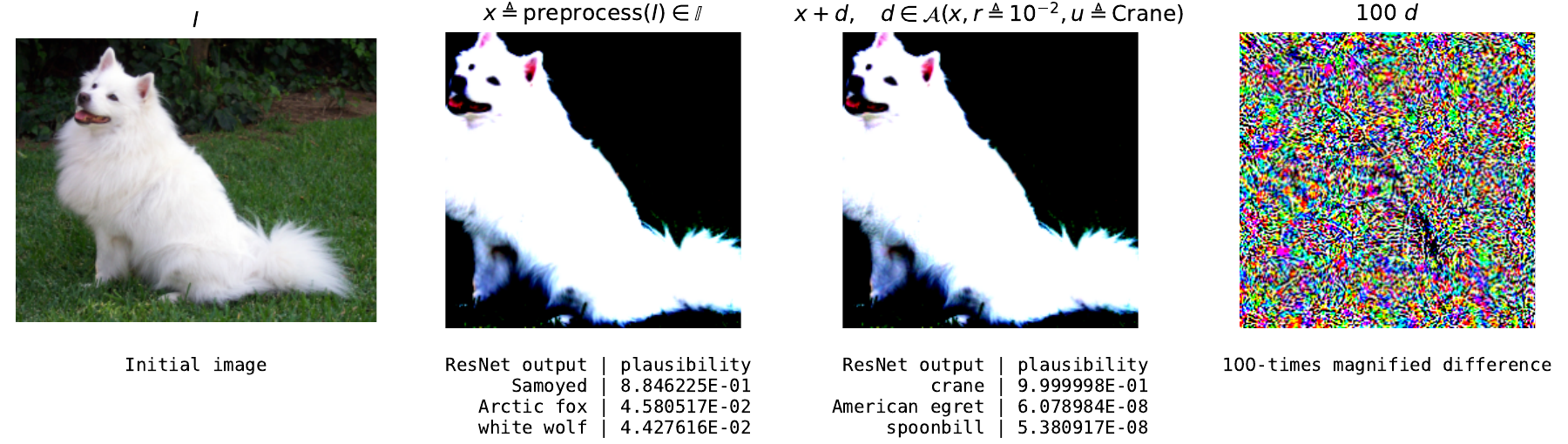}
    \caption{
        Targeted attack on \resnet.
        (Left) Image of a Samoyed dog.
        (Centre left) Preprocessed image and its classification.
        (Centre right) Attack of the preprocessed image, targeting the class~$"\mathrm{Crane}"$ and allowing to alter each pixel by at most~$10^{-2}$ units, and its classification.
        (Right) Magnification of the image alteration performed by the attack.
        }
    \label{figure:illustration_attack}
\end{figure}

Consider PyTorch's \href{https://pytorch.org/hub/pytorch_vision_resnet/}{\resnet} \nn for image classification (with its default training).
This network classifies images according to~$1000$ pre-selected classes.
That is,~$\resnet$ takes as inputs RGB images with~$224 \times 224$ pixels, represented as tensors in~$\Ibb \defequal [0,1]^{3 \times 224 \times 224}$, and outputs vectors in the 1000-dimensional space~$\Pbb^{1000} \defequal \{p \in [0,1]^{1000}: \sum_{\ell=1}^{1000} p_\ell = 1\}$, where each component represents the plausibility that the image depicts an item from the associated class.
Then the image is classified to the class~$\ell \in \llb1,1000\rrb$ with the highest value~$p_\ell$.
Consider the image in Figure~\ref{figure:illustration_attack} (left), which depicts a Samoyed dog.
After preprocessing, the image becomes a tensor~$x \in \Ibb$ (centre left), and \resnet correctly classifies it as a Samoyed with~$88\%$ confidence.
Then, consider a targeted attack shifting the classification of~$x$ from a Samoyed to a crane.
This consists in a perturbation~$d$, with~$\textnorm{d}_\infty \leq 10^{-2}$, that shifts~$\resnet(x+d)$ towards the one-hot vector associated to the~$"\mathrm{Crane}"$ class.
The altered image~$x+d$ (centre right) remains visually indistinguishable from the original, yet the classification changes drastically: \resnet assigns nearly~$100\%$ confidence to the \inlinequote{$\mathrm{Crane}$} class.
As shown on the right of the figure,~$d$ alters many pixels, but by at most~$10^{-2}$ units so the result is visually unchanged.

\subsubsection{Formal Definition}
\label{section:optimization_attack/notion/definition}

We now formalize directional \nn attacks in Definition~\ref{definition:directional_attack}.
We build them from targeted \nn attacks, which we therefore introduce first.
For simplicity, we state both notions directly in terms of the \nn~$\Phi$ and the constraints function~$c$.
Their formulations also involve a loss function~$\Lcal: \Rbb^m \times \Rbb^m \to \Rbb_+$ and a norm~$\textnorm{\cdot}$, that are left generic since several usual choices (that we discuss in Remark~\ref{remark:usual_choices_attack}) may be  considered.
We also adapt Definition~\ref{definition:directional_attack} to \emph{directional \nn attacks with respect to subspaces} in Remark~\ref{remark:directional_attack_components}.

The notion of targeted \nn attack must be adapted from the classification contexts.
Consider a \emph{loss} function~$\Lcal: \Rbb^m \times \Rbb^m \to \Rbb_+$ and a norm~$\textnorm{\cdot}$.
For all~$(x,r,y) \in \Rbb^n \times \Rbb_+ \times \Rbb^m$, define the problem
\begin{equation*}
    \label{problem:targeted_attack}
    \problemoptimoneline
        {\minimize}
        {d \in \Bcal_r^n}
        {\Lcal\left(\Phi(x+d) , y\right)}
        {c(\Phi(x+d)) \leq 0.}
\end{equation*}
A \emph{targeted attack on~$\Phi$ at~$x$ with radius~$r$ and target~$y$} (with respect to~$\Lcal$ and the ball~$\Bcal_r^n$ in the~$\textnorm{\cdot}$ norm) consists in an input alteration~$d \in \Bcal_r^n$ chosen with respect to this problem.
All feasible elements are said to be \emph{feasible attacks}, all global solutions are said to be \emph{optimal attacks}, and all feasible attacks~$d$ satisfying~$\Lcal(\Phi(x+d),y) < \Lcal(\Phi(x),y)$ are said to be \emph{successful attacks}.
From that notion, we may formalize those of a \emph{directional \nn attack} as Definition~\ref{definition:directional_attack}.

\begin{definition}[Directional \nn attack]
    \label{definition:directional_attack}
    Consider a \emph{loss} function~$\Lcal: \Rbb^m \times \Rbb^m \to \Rbb_+$ and a norm~$\textnorm{\cdot}$.
    For all~$x \in \Rbb^n$, define the \emph{differential \nn}~$\Phi_x: d \in \Rbb^n \mapsto \Phi(x+d)-\Phi(x) \in \Rbb^m$.
    Then, for all points~$x \in \Rbb^n$, all radii~$r \in \Rbb_+$, and all directions~$u \in \Rbb^m$, an \emph{attack of~$\Phi$ at~$x$ of radius~$r$ in direction~$u$} (with respect to~$\Lcal$ and the ball~$\Bcal_r^n$ in the~$\textnorm{\cdot}$ norm) consists in a targeted attack on~$\Phi_x$ at~$0$ with radius~$r$ and target~$u\textnorm{u}^{-1}$ if~$u \neq 0$ and~$0$ otherwise.
    That is, directional attacks are input alterations~$d \in \Bcal_r^n$ chosen with respect to the problem \inlinequote{$\minimize_{d \in \Bcal_r^n} \Lcal(\Phi_x(d),0)$ subject to~$c(\Phi(x)+\Phi_x(d)) \leq 0$} if~$u = 0$, and to the problem
    \begin{equation}
        \label{problem:directional_attack}
        \tag{$\Abf_\Lcal(x, r, u)$}
        \problemoptimoneline
            {\minimize}
            {d \in \Bcal_r^n}
            {\Lcal\left(\Phi_x(d) , \dfrac{u}{\norm{u}}\right)}
            {c(\Phi(x)+\Phi_x(d)) \leq 0}
    \end{equation}
    if~$u \neq 0$.
    We denote by~$\attacksetfeas(x, r, u)$ the set of feasible attacks, by~$\attacksetsuccess(x, r, u)$ the set of successful attacks, and by~$\attacksetopt(x, r, u)$ the set of optimal attacks.
\end{definition}

\begin{remark}
    \label{remark:usual_choices_attack}
    As Definition~\ref{definition:directional_attack} stresses, all solvers from the literature designed for targeted attacks may be used off-the-shelf to compute directional attacks, since the latter is a specific instance of the former.
    Moreover, Definition~\ref{definition:directional_attack} is compatible with any loss function and any norm.
    However, many numerical tools are restricted to the~$\textnorm{\cdot}_\infty$ norm and, as shown in Section~\ref{section:optimization_attack/notion/numerical_tools}, either the \emph{square-error loss} function~$\squarerrorloss$ or the \emph{cross-entropy loss} function~$\crossentropyloss$.
    These two losses are defined as follows, for all~$(y_1,y_2) \in \Rbb^m \times \Rbb^m$, and denoting by~$\vv{\ln}$ the logarithm applied component-wise,
    \begin{equation*}
        \squarerrorloss(y_1,y_2) \defequal \norm{y_2-y_1}_2^2
        \quad \mbox{and} \quad
        \crossentropyloss(y_1,y_2) \defequal -\dotprod{\vv{\ln}(\softmax(y_1))}{\softmax(y_2)}.
    \end{equation*}
\end{remark}

\begin{remark}
    \label{remark:directional_attack_components}
    Definition~\ref{definition:directional_attack} is designed so that, at any~$x$ and any direction~$u$, an optimal attack direction~$d$ makes the vector~$\Phi_x(d)$ to match as much as possible with~$u$ (according to the chosen loss function~$\Lcal$).
    However, in some contexts, we may want for~$\Phi_x(d)$ to match~$u$ with respect to some components only, or more generally with respect to their projections to some linear subspace.
    Given a~$(m \times m)$ orthogonal projection matrix~$P$, we seek a direction~$d \in \Rbb^n$ such that the projection~$P\Phi_x(d)$ agrees with the projection~$Pu$.
    A \emph{directional \nn attack of~$\Phi$ at~$x$ in the direction~$u$ with respect to~$P$} then consists in attacking (in the sense of Definition~\ref{definition:directional_attack}) the \emph{projected differential \nn~$P\Phi_x$} in the direction~$Pu$.
\end{remark}

\subsubsection{Numerical Tools}
\label{section:optimization_attack/notion/numerical_tools}

Several numerical tools are available for computing \nn attacks.
Most of them are implemented in PyTorch~\cite{PyTorch}, and involve one of the loss functions defined in Remark~\ref{remark:usual_choices_attack}.
The open-source solvers we are aware of are listed in Table~\ref{table:list_solvers_nn_attacks}, filled to our best knowledge by investigating their codes and documentations.
While these tools are designed for targeted attacks only, Definition~\ref{definition:directional_attack} shows how to use them to compute directional attacks.
In addition, all solvers except~\abcrown are compatible with \nn models not given as glass-boxes, as they only require that the \nn allows for backpropagation (none of the solvers implement algorithms~\cite{AlShChZhHsSr19GenAttack,AnCrFlHe20SquareAttack,GuGaYoWiWe19SimpleBlackBoxAttack} for attacks of black-box \nns).
However, none of these solvers currently supports constraints on the input, so we introduce a workaround in Section~\ref{section:optimization_attack/use_in_optim/heuristic_tool} to address this limitation.

\begin{table}[!ht]
    \centering
    \begin{tabular}{| R{0.18\linewidth} | M{0.08\linewidth} | L{0.60\linewidth} |}
        \hline
        Software & Losses & Additional requirements \\ \hline \hline
        \abcrown~\cite{zhang22babattack} & $\squarerrorloss$ & $c \equiv 0$,~$\Phi$ glass-box \relu-based model \\ \hline
        FoolBox~\cite{rauber2017foolbox} & $\crossentropyloss$ & $c \equiv 0$, backpropagation \\ \hline
        ART~\cite{AdvRobustToolbox} & $\crossentropyloss$ & $c \equiv 0$, backpropagation \\ \hline
        AdvertTorch~\cite{AdvertTorch} & $\crossentropyloss$ & $c \equiv 0$, backpropagation \\ \hline
        CleverHans~\cite{CleverHans} & $\crossentropyloss$ & $c \equiv 0$, backpropagation \\ \hline
        Torchattacks~\cite{kim2020torchattacks} & $\crossentropyloss$ & $c \equiv 0$, backpropagation,~$\Phi: [0,1]^n \to \Rbb^m$ \\ \hline
    \end{tabular}
    \caption{PyTorch-compatible solvers for directional \nn attacks we are aware of.}
    \label{table:list_solvers_nn_attacks}
\end{table}

\subsection{Optimizing Through Directional \nn Attacks}
\label{section:optimization_attack/use_in_optim}

We now leverage directional \nn attacks to solve Problem~\eqref{problem:P}.
For all points~$x \in \Rbb^n$, a well-designed directional attack of~$\Phi$ at~$x$ in the direction~$\nabla f(\Phi(x))$ yields an ascent direction for Problem~\eqref{problem:P}.
We introduce the \emph{\attackstep operator}, which performs the aforementioned attack, in Definition~\ref{definition:attack_operator}.
The norm is left unspecified since it does not impact our analysis, but in practice we consider the~$\textnorm{\cdot}_\infty$ norm.

\begin{definition}[\attackstep operator]
    \label{definition:attack_operator}
    A set-valued function~$\attackstep: \Rbb^n \times \Rbb_+ \to 2^{\Rbb^n}$ is an \emph{attack operator} for some loss function~$\Lcal$ if for all~$(x,r) \in \Rbb^n \times \Rbb_+$, it holds that~$\attackstep(x,r) \subseteq \attacksetfeas(x, r, \nabla f(\Phi(x)))$.
\end{definition}

Many functions satisfy Definition~\ref{definition:attack_operator}, e.g., the empty-set operator~$\attackstep \equiv \emptyset$ and the ideal operator given, for all~$(x,r) \in \Rbb^n \times \Rbb_+$, by~$\attackstep(x,r) \defequal \attacksetopt(x, r, \nabla f(\Phi(x)))$.
In practice, we choose a numerical solver from the literature, and for all~$(x,r) \in \Rbb^n \times \Rbb_+$ we define~$\attackstep(x,r)$ as the output of that solver tackling Problem~\eqref{problem:directional_attack} at~$u \defequal \nabla f(\Phi(x))$.
The solver either fails (so it returns~$\emptyset$, which is allowed) or returns elements of~$\attacksetfeas(x,r,u)$ that may not be optimal or even successful.
Hence, we do not impose that~$\attackstep(x,r) \subseteq \attacksetopt(x,r,u)$ or~$\attackstep(x,r) \subseteq \attacksetsuccess(x,r,u)$.

The~$\attackstep$ operator acts as a suitable replacement of the chain rule to calculate gradients of~$f \circ \Phi$, even when the Jacobian of~$\Phi$ is nonexistent or in the presence of constraints.
It also generalizes the computation of ascent directions by backpropagations of gradients of~$f$ through~$\Phi$, in the sense that
\begin{equation*}
    c \equiv 0
    \implies
    \fct{\attackstep}
        {(x,r)}
        {\Rbb^n \times \Rbb_+}
        { \left\{ \dfrac{r \ \backprop(\Phi,x,\nabla f(\Phi(x)))}{\norm{\backprop(\Phi,x,\nabla f(\Phi(x)))}} \right\} }
        {2^{\Rbb^n}}
    ~\text{satisfies Definition~\ref{definition:attack_operator} for}~ \squarerrorloss.
\end{equation*}
Moreover, even in unconstrained cases, calculating ascent directions from directional attack algorithms may outperform direct backpropagations of~$\nabla f(\Phi(x))$, especially in cases where~$\Phi$ is prone to numerical difficulties such as the \textit{vanishing gradient}~\cite{Ho98VanishingGradient} or when~$\textnorm{\nabla f(\Phi(x))}$ becomes small.
We substantiate this claim in Section~\ref{section:numerical}.

In Section~\ref{section:optimization_attack/use_in_optim/setting_ascent_directions}, we prove that if the~$\attackstep$ operator is actually guaranteed to identify successful attacks, then its outputs possess guarantees to be ascent directions for Problem~\eqref{problem:P}.
In Section~\ref{section:optimization_attack/use_in_optim/heuristic_tool}, we discuss that although this setting is not fully supported by existing solvers, an \attackstep operator defined via current \nn attacks solvers already acts as a good heuristic for computing ascent directions.

\subsubsection{Favourable Setting for the \attackstep Operator}
\label{section:optimization_attack/use_in_optim/setting_ascent_directions}

This section establishes a theoretical setting guaranteeing that the \attackstep operator yields an ascent direction for Problem~\eqref{problem:P} whenever one exists.
This setting requires that the~$\attackstep$ operator returns successful attack directions when some exist, and that the loss function is \emph{well-suited} in a sense that we introduce below.
We formalize these requirements in Assumption~\ref{assumption:attack_guarantees_ascent_directions} and our claim in Proposition~\ref{proposition:careful_attack_yield_dominating_direction}.

We say that a function~$\Lcal: \Rbb^m \times \Rbb^m \to \Rbb_+$ satisfies the \emph{Well-Suited Loss Property} (\wslp) when
\begin{equation}
    \label{equation:well_suited_losses}
    \tag{\wslp}
    \exists c > 0: \quad
    \forall (y_1,y_2) \in \Rbb^m \times \Rbb^m, \quad
    \Lcal(y_1,y_2) < \Lcal(0,y_2) \implies \textdotprod{y_1}{y_2} > c\textnorm{y_1}.
\end{equation}
Among the loss functions from Remark~\ref{remark:usual_choices_attack},~$\squarerrorloss$ satisfies the \wslp while~$\crossentropyloss$ does not (see Remark~\ref{remark:usual_losses_wslp}).
Roughly speaking, the \wslp ensures that for all~$(x,r) \in F \times \Rbb_+$, all attacks~$d \in \attackset^+(x, r, \nabla f(\Phi(x)))$ are alterations of~$x$ that drive the output of~$\Phi$ in the direction of the gradient of~$f$.
Thus, the first-order approximation of~$f$ near~$\Phi(x)$ ensures that if~$\textnorm{d}$ is small enough, then~$f(\Phi(x+d)) > f(\Phi(x))$.

\begin{assumption}[Successful \attackstep operator]
    \label{assumption:attack_guarantees_ascent_directions}
    The loss function~$\Lcal$ defining the~$\attackstep$ operator satisfies the \wslp, and for all~$(x,r) \in \Rbb^n \times \Rbb_+$, it holds that~$\attackstep(x,r) \subseteq \attacksetsuccess(x,r,\nabla f(\Phi(x)))$.
    Moreover, for all~$(x,r) \in \Rbb^n \times \Rbb_+$ such that~$\attacksetsuccess(x,r,\nabla f(\Phi(x)))$ is nonempty,~$\attackstep(x,r)$ is nonempty.
\end{assumption}

\begin{proposition}
    \label{proposition:careful_attack_yield_dominating_direction}
    Under Assumptions~\ref{assumption:problem} and~\ref{assumption:attack_guarantees_ascent_directions}, for all~$x \in \Rbb^n$, there exists~$\overline{r}(x) > 0$ such that for all~$r \in [0,\overline{r}(x)]$, every~$d \in \attackstep(x, r)$ is an ascent direction for Problem~\eqref{problem:P} emanating from~$x$.
\end{proposition}

\begin{proof}
    Let Assumptions~\ref{assumption:problem} and~\ref{assumption:attack_guarantees_ascent_directions} hold.
    Let~$c > 0$ be the constant involved in the \wslp, and define
    \begin{equation*}\begin{array}{rlcl}
        \forall y \in \Rbb^m,
            & \overline{\rho}(y)
            & \defequal
            & \max\{ \rho \in [0,+\infty] : \quad (\forall u \in \Bcal^m_\rho,~ \textdotprod{u}{\nabla f(y)} > c\textnorm{u} \implies f(y+u) > f(y)) \}, \\
        \forall x \in \Rbb^n,
        & \overline{r}(x)
        & \defequal
        & \max\{ r \in [0,+\infty] : \quad (\forall d \in \Bcal^n_r,~ \textnorm{\Phi_x(d)} \leq \overline{\rho}(\Phi(x))) \}.
    \end{array}\end{equation*}
    First, for all~$y \in \Rbb^m$, the first-order approximation of~$f$ near~$y$ ensures that~$\overline{\rho}(y) > 0$.
    Moreover, for all~$x \in \Rbb^n$, the continuity of~$\Phi$ from Assumption~\ref{assumption:problem} ensures that~$\overline{r}(x) > 0$.
    We deduce that
    \begin{equation}
        \label{equation:property_in_proof_careful_attacks}
        \tag{$\star$}
        \forall d \in \Bcal^n_{\overline{r}(x)}: \dotprod{\Phi_x(d)}{\nabla f(\Phi(x))} > c\norm{\Phi_x(d)}, \qquad
        f(\Phi(x+d)) > f(\Phi(x)).
    \end{equation}
    If moreover~$c(\Phi(x+d)) \leq 0$, then~$d$ is an ascent direction emanating from~$x$.
    Second, Assumption~\ref{assumption:attack_guarantees_ascent_directions} ensures that for all~$x \in \Rbb^n$ and all~$r \in [0,\overline{r}(x)]$, all~$d \in \attackstep(x,r) \subseteq \attacksetsuccess(x,r,\nabla f(\Phi(x)))$ satisfy
    \begin{equation*}
        \begin{array}{rrl}
            & c(\Phi(x+d)) \leq 0 & \text{since}~ d ~\text{is feasible}, \\[0.5em]
            & d \in \Bcal^n_{\overline{r}(x)} & \text{since}~ \norm{d} \leq r \leq \overline{r}(x), \\[0.5em]
            & \Lcal(\Phi_x(d), \nabla f(\Phi(x))) < \Lcal(0, \nabla f(\Phi(x))) & \text{since}~ d ~\text{is a successful attack},
            \\
            \mbox{so}
            & \dotprod{\Phi_x(d)}{\nabla f(\Phi(x))} > c\norm{\Phi_x(d)} & \text{since}~ \Lcal ~\mbox{satisfies the \ref{equation:well_suited_losses}}.
        \end{array}
    \end{equation*}
    The result follows directly thanks to~\eqref{equation:property_in_proof_careful_attacks} and the guarantee that~$c(\Phi(x+d)) \leq 0$.
\end{proof}

Proposition~\ref{proposition:careful_attack_yield_dominating_direction} validates the first goal of this paper: directional \nn attacks may be used as tools for solving Problem~\eqref{problem:P}.
Under Assumption~\ref{assumption:attack_guarantees_ascent_directions}, for all~$x \in \Rbb^n$ and all~$r \leq \overline{r}(x)$, all~$d \in \attackstep(x,r)$ are ascent directions for Problem~\eqref{problem:P}.
However, note that, at the time of writing, no \attackstep operator constructed from a solver in Table~\ref{table:list_solvers_nn_attacks} satisfies Assumption~\ref{assumption:attack_guarantees_ascent_directions} (we discuss this in Section~\ref{section:optimization_attack/use_in_optim/heuristic_tool}).
Yet, as we show in Section~\ref{section:optimization_hybrid}, the \attackstep operator remains valuable within a broader optimization strategy.

We conclude this section with several remarks.
Remark~\ref{remark:maximal_radius_attack_guarantees_ascent_directionss} shows how to search for a radius~$r \leq \overline{r}(x)$, Remark~\ref{remark:empty_set_of_attacks} discusses whether~$\attackstep(x,r) = \emptyset$, and Remark~\ref{remark:which_ascent_direction_to_attack} adapts our methodology to cases where~$f$ has a Hessian and cases where~$f$ has no gradient.
Finally, Remark~\ref{remark:active_subspaces} focuses on cases where~$f$ has \emph{active subspaces}, which require care to avoid a deterioration of the performance of the \attackstep operator.

\begin{remark}
    \label{remark:maximal_radius_attack_guarantees_ascent_directionss}
    Under Assumption~\ref{assumption:attack_guarantees_ascent_directions}, Proposition~\ref{proposition:careful_attack_yield_dominating_direction} shows that, given any~$x \in \Rbb^n$, we must select~$r \leq \overline{r}(x)$ to ensure that all~$d \in \attackstep(x,r)$ are ascent directions.
    Computing~$\overline{r}(x)$ is difficult, since it depends on the local Lipschitz constant of~$\Phi$ at~$x$.
    However, it is strictly positive, so we may find some~$r \in {]}0,\overline{r}(x)]$ by initializing~$r \defequal 1$ and halving it until any~$d \in \attackstep(x,r)$ is an ascent direction.
\end{remark}

\begin{remark}
    \label{remark:empty_set_of_attacks}
    Some~$x \in \Rbb^n$ yield~$\attackstep(x,r) = \emptyset$ for all~$r \geq 0$ small enough, so Proposition~\ref{proposition:careful_attack_yield_dominating_direction} is vacuously true at such points.
    This happens, for example, when~$(x,r)$ satisfies a necessary optimality conditions for Problem~\eqref{problem:P} given by~$\textdotprod{\Phi_x(d)}{\nabla f(\Phi(x))} \leq 0$ for all~$d \in \Bcal_r^n$ such that~$x+d \in F$.
    It may also happen that no successful attack exists when~$(x,r)$ does not satisfy this necessary optimality condition (so~$\attackstep(x,r) = \emptyset$).
    Consider the following example, where~$\mathrm{Id}_{\Rbb^2}$ denotes the identity function in~$\Rbb^2$,
    \begin{equation*}
        \Lcal \defequal \squarerrorloss, \quad
        f \defequal (y \in \Rbb^2 \mapsto y_2), \quad
        c \defequal \left(y \in \Rbb^2 \mapsto y_2-\dfrac{y_1^2}{2}\right), \quad
        \Phi \defequal \mathrm{Id}_{\Rbb^2}, \quad
        x \defequal (0,0), \quad
        r \defequal 1.
    \end{equation*}
    In this context, it holds that~$\Lcal(\Phi_x(d), \nabla f(\Phi(x))) < \Lcal(\Phi_x(0), \nabla f(\Phi(x)))$ if and only if~$\textnorm{d-(0,1)} < 1$.
    However, no such direction~$d$ yields~$x+d$ in the feasible set~$F = \{x \in \Rbb^2: x_2 \leq \frac{x_1^2}{2}\}$.
    As a consequence, no feasible attack exists.
    Yet, there exists ascent directions for~$f \circ \Phi$ in~$\Bcal_r(x)$, such as~$d \defequal (\frac{1}{2}, \frac{1}{8})$.
\end{remark}

\begin{remark}
    \label{remark:which_ascent_direction_to_attack}
    The~$\attackstep$ operator computes directional \nn attacks at all~$x \in \Rbb^n$ in the \emph{gradient ascent} direction~$\nabla f(\Phi(x))$ only.
    Yet, other directions may be considered depending on the regularity of~$f$, such as the \emph{Hessian ascent direction}~$[\nabla^2 f(\Phi(x))]^{-1} \nabla f(\Phi(x))$ when~$f$ has a Hessian~$\nabla^2 f$.
    When~$f$ is nonsmooth, \emph{subgradient}, \emph{simplex gradient}~\cite{HaJBPl23PositiveBases,HaJBPl20ErrorsBoundsSimplexGrad} or even \emph{simplex Hessian}~\cite{HaJBCh20Hessian} may be considered for ascent directions.
    For all~$y \in \Rbb^m$ and all~$r \in \Rbb_+^*$, the (canonical forward) simplex gradient of~$f$ of radius~$r$ at~$y$ is the vector~$\nabla_r f(y) \defequal r^{-1}[f(y+re_i)-f(y)]_{i=1}^{m}$, and the (canonical forward) simplex Hessian of~$f$ of radius~$r$ at~$y$ is the matrix~$\nabla_r^2 f(y) \defequal r^{-1}[\nabla_r f(y+re_j) - \nabla_r f(y)]_{j=1}^{m}$.
\end{remark}

\begin{remark}
    \label{remark:active_subspaces}
    The \attackstep operator falls short of its best potential when~$f$ has an active subspace~\cite{active_subspace_perf,active_subspaces} of dimension~$p \ll m$.
    We say that~$f$ has an active subspace if there exists a mapping~$g: \Rbb^m \to \Rbb$ and an orthogonal projection matrix~$P \in \Rbb^{m \times m}$ with rank~$p$ such that~$f(y) = g(Py)$ for all~$y \in \Rbb^m$.
    In this case, for all~$x \in \Rbb^n$, we have~$\nabla f(\Phi(x)) = P \nabla g(P\Phi(x))$, so a directional \nn attack from Definition~\ref{definition:directional_attack} seeks for a direction~$d \in \Rbb^n$ such that~$\Phi_x(d) \approx P\nabla g(P\Phi(x))$.
    However, the components of~$\Phi_x(d)$ orthogonal to the span of~$P$ are actually irrelevant, since~$f(\Phi(x)+\Phi_x(d)) = f(\Phi(x) + P\Phi_x(d))$.
    If~$P$ is known, the \attackstep operator can be adapted as in Remark~\ref{remark:directional_attack_components}.
    Otherwise, we may learn~$P$ on the fly~\cite{active_subspaces_learning} and adapt accordingly.
    Active subspaces arise naturally, for example, when~$\Phi$ has the form~$\Phi(x) = [\Psi(x), x]$ for all~$x \in \Rbb^n$, with~$\Psi: \Rbb^n \to \Rbb^k$ another \nn, allowing to model constraints on~$x$ via the function~$c$.
    In such cases, the objective may depend only on~$\Psi(x)$, so the function~$f: \Rbb^{k+n} \to \Rbb$ has an active subspace of dimension~$k$ (with~$P$ the diagonal matrix with~$k$ ones followed by~$n$ zeros).
\end{remark}

\subsubsection{Practical Construction of \attackstep Operators}
\label{section:optimization_attack/use_in_optim/heuristic_tool}

Proposition~\ref{proposition:careful_attack_yield_dominating_direction} describes an \attackstep operator that returns ascent directions for Problem~\eqref{problem:P}.
However, its setting (enclosed by Assumption~\ref{assumption:attack_guarantees_ascent_directions}) may not hold in practice, where we construct an \attackstep operator by selecting a solver~$\Srm$ and, for all~$(x,r) \in \Rbb^n \times \Rbb_+$, defining~$\attackstep(x,r)$ as the output of~$\Srm$ solving Problem~\eqref{problem:directional_attack} at~$u \defequal \nabla f(\Phi(x))$.
This section discusses the gap between Assumption~\ref{assumption:attack_guarantees_ascent_directions} and practical settings, and introduces a heuristic yet efficient \attackstep operator from existing solvers.
This heuristic operator is stated in Definition~\ref{definition:attack_operator_practical}, and its performance is assessed in Section~\ref{section:numerical}.

Assumption~\ref{assumption:attack_guarantees_ascent_directions} is currently incompatible with solvers listed in Table~\ref{table:list_solvers_nn_attacks}, for four independent reasons.
Below, we introduce a workaround to three of them.
The only incompatibility we cannot address is the requirement that the \attackstep operator returns a successful attack whenever one exists.
To our best knowledge, only \abcrown~\cite{zhang22babattack} ensures this (but under the requirement that~$\Phi$ is a glass-box).

The first problem we encountered is that no available solver supports input constraints.
Our workaround to address this drawback consists in a (heuristic) reformulation of the attack problem based on quadratic penalization.
Given a penalization weight~$\lambda$ (with~$\lambda \defequal 10$ by default), we define
\begin{equation*}
    \fct
        {\widetilde{\Phi}}
        {x}
        {\Rbb^n}
        {\begin{bmatrix} \Phi(x) \\ \relu(c(\Phi(x))) \end{bmatrix}}
        {\Rbb^{m+p}}
    \qquad \mbox{and} \qquad
    \fct
        {\widetilde{f}}
        {\begin{bmatrix} y \\ z \end{bmatrix}}
        {\Rbb^{m+p}}
        {f(y)-\lambda\norm{z}_2^2.}
        {\Rbb}
\end{equation*}
Then, for any~$(x, r) \in \Rbb^n \times \Rbb_+$, we do not compute directional attacks on~$\Phi$ at~$x$ of radius~$r$ in direction~$\nabla f(\Phi(x))$ by solving Problem~\eqref{problem:directional_attack} at~$u \defequal \nabla f(\Phi(x))$ (with constraints given explicitly in the attack problem).
Instead, we solve the (unconstrained) problem of minimizing~$\Lcal(\widetilde{\Phi}_x(d), \nabla \widetilde{f}(\widetilde{\Phi}(x)))$ over~$d \in \Bcal_r^n$, where~$\widetilde{\Phi}_x$ is the differential \nn of~$\widetilde{\Phi}$ defined by~$\widetilde{\Phi}_x(d) \defequal \widetilde{\Phi}(x+d)-\widetilde{\Phi}(x)$ for all~$d \in \Rbb^n$.
This problem is tractable using current solvers since it is unconstrained, but the solution it returns may be infeasible.

The second incompatibility is that some solvers assume input domains restricted to~$[0,1]^n$ rather than~$\Rbb^n$.
In this case, we define the affine function~$d_r: \delta \mapsto 2r(\delta - \1/2)$ that maps~$[0,1]^n$ to the ball~$\Bcal^n_r$ in the~$\textnorm{\cdot}_\infty$ norm.
So, for all~$(x,r) \in \Rbb^n \times \Rbb_+$ we define the scaled differential \nn
\begin{equation*}
    \fct
        {\Phi_{(x,r)}}
        {\delta}
        {[0,1]^n}
        {\Phi_x(d_r(\delta)),}
        {\Rbb^m}
\end{equation*}
so, for any~$(x, r, u) \in \Rbb^n \times \Rbb_+ \times \Rbb^m$, we can get a directional attack on~$\Phi$ at~$x$ of radius~$r$ in direction~$u$ by computing a targeted attack on~$\Phi_{(x,r)}$ at~$\1/2$ with radius~$1/2$ and target~$u$.
This new problem is compatible with the input domain restrictions of all solvers we listed in Table~\ref{table:list_solvers_nn_attacks}, and any direction~$\delta$ has the same classification (optimal, successful, or unsuccessful) for this problem as the direction~$d_r(\1/2+\delta) = 2r\delta$ for the original directional attack of~$\Phi$ at~$x$ with radius~$r$.

Finally, the third drawback of most solvers is their reliance on the cross-entropy loss function~$\crossentropyloss$, which does not satisfy the \wslp.
However, we have found that the composition of~$\crossentropyloss$ with a some projection remains principled with respect to the \wslp.
Indeed (we refer to Remark~\ref{remark:usual_losses_wslp} for details),
\begin{equation*}
    \forall (y_1, y_2) \in \Rbb^m \times \Rbb^m \setminus \{0\}, \quad
    \crossentropyloss\left(P(y_2)y_1, P(y_2)y_2\right) < \crossentropyloss\left(0, P(y_2)y_2\right) \implies \textdotprod{y_1}{y_2} > 0,
\end{equation*}
where for all~$y \in \Rbb^m \setminus \{0\}$, we define the projection matrix~$P(y) \defequal y\transpose{y} \textnorm{y}^{-2}$.
This logical implication does not match the \wslp, but it is principled in the sense that the formula claims that~$\textdotprod{y_1}{y_2} > 0$ instead of~$\textdotprod{y_1}{y_2} > c\textnorm{y_1}$.
Then, we may use off-the-shelf the solvers relying on~$\crossentropyloss$ by defining the attack problem as in Remark~\ref{remark:directional_attack_components}.
For all~$(x,u) \in \Rbb^n \times \Rbb^m \setminus \{0\}$, we may compute attacks of~$\Phi$ at~$x$ in the direction~$u$ by using any solver to compute targeted attacks of~$P(u)\Phi_x$ at~$0$ with target~$P(u)u$.

We now build an~$\attackstep$ operator from these workarounds, in Definition~\ref{definition:attack_operator_practical}.
It is heuristic when Problem~\eqref{problem:P} has constraints, since it then works on inexact relaxation of Problem~\eqref{problem:P}.

\begin{definition}[practical \attackstep operator]
    \label{definition:attack_operator_practical}
    Let~$\Srm$ be a solver for targeted \nn attacks and let~$\lambda \geq 0$.
    The \emph{$\attackstep$ operator induced by~$\Srm$} is the set-valued function~$\attackstep_{(\Srm,\lambda)}: \Rbb^n \times \Rbb_+ \to 2^{\Rbb^n}$ that, for all~$(x,r) \in \Rbb^n \times \Rbb_+$, maps~$(x,r)$ to the outputs of~$\Srm$ solving the problem of attacking the \nn~$\smash{\widetilde{P} \widetilde{\Phi}_{(x,r)}}$ at~$\1/2$ with radius~$1/2$ and target~$\smash{\widetilde{P} \nabla \widetilde{f}(\widetilde{\Phi}(x))}$, where~$\smash{\widetilde{P} \defequal P(\nabla \widetilde{f}(\widetilde{\Phi}(x)))}$, using the loss~$\Lcal_\Srm$ fixed by~$\Srm$.
    That is,~$\attackstep_{(\Srm,\lambda)}(x,r)$ returns the singleton~$\attackstep_{(\Srm,\lambda)}(x,r) \defequal \{2r\delta\}$, where~$\delta$ is the output of~$\Srm$ solving the problem
    \begin{equation*}
        \problemoptimfree
            {\minimize}
            {\delta \in \Bcal_{1/2}^n}
            {\Lcal_\Srm\left( \widetilde{P} \widetilde{\Phi}_{(x,r)}(\1/2+\delta) , \widetilde{P} \nabla\widetilde{f}(\widetilde{\Phi}(x))\right).}
    \end{equation*}
\end{definition}

We emphasize that the current gap between theoretical and practical \attackstep operators may narrow as solver implementations advance.
In our experiments in Section~\ref{section:numerical}, we define the~$\attackstep_{(\Srm,\lambda)}$ operator with the Torchattacks library~\cite{kim2020torchattacks} as the solver~$\Srm$ (specifically, its implementations of the \fgsm~\cite{Goodfellow15FGSM} and \pgd~\cite{Madry19PGD} algorithms under the~$\textnorm{\cdot}_\infty$ norm and the~$\crossentropyloss$ loss function).

\begin{remark}
    \label{remark:usual_losses_wslp}
    The~$\squarerrorloss$ loss function satisfies the \wslp, with~$c \defequal \frac{1}{2}$.
    Indeed,
    \begin{equation*}
        \begin{array}{rrcl}
            \forall (y_1, y_2) \in \Rbb^m \times \Rbb^m, \quad
                & \squarerrorloss(y_1,y_2) < \squarerrorloss(0,y_2)
                & \iff       & \textdotprod{y_2-y_1}{y_2-y_1} < \textdotprod{y_2}{y_2} \\
                & & \iff     & -2\textdotprod{y_1}{y_2} + \textdotprod{y_1}{y_1} < 0 \\
                & & \implies & \textdotprod{y_1}{y_2} > \frac{1}{2}\textnorm{y_1}.
        \end{array}
    \end{equation*}
    In contrast,~$\crossentropyloss$ admits counter-examples to the \wslp, such as
    \begin{equation*}
        \begin{cases}
            y_1 \defequal \left(-0.1, -1\right), \\
            y_2 \defequal \left(   1,  0\right),
        \end{cases}
        \quad \text{which yields} \quad
        \begin{cases}
            \text{either} & \crossentropyloss(y_1, y_2) \approx 0.58 < \crossentropyloss(0, y_2) \approx 0.69 \\
            \text{and} & \textdotprod{y_1}{y_2} = -0.1 \leq 0.
        \end{cases}
    \end{equation*}
    Yet, an alteration of~$\crossentropyloss$ shares similarities with the \wslp.
    Consider, for all~$(y_1,y_2) \in \Rbb^m \times \Rbb^m \setminus \{0\}$, the loss~$\crossentropyloss(P(y_2)y_1, P(y_2)y_2)$, where we define the projection~$P(y) \defequal y\transpose{y}\textnorm{y}^{-2}$.
    Denote, for all~$y \in \Rbb^m$, by~$\mathrm{le}(y) \defequal \ln(\sum_{i=1}^{m} e^{y_i})$, and remark that, for all~$(y_1, y_2) \in \Rbb^m \times \Rbb^m \setminus \{0\}$, we may rewrite~$\crossentropyloss(y_1,y_2)$ as~$\crossentropyloss(y_1, y_2) = \mathrm{le}(y_1) - \textdotprod{y_1}{\softmax(y_2)}$ and~$P(y_2)y_1 = \textnorm{y_2}^{-2}\textdotprod{y_1}{y_2}y_2$ and~$P(y_2)y_2 = y_2$.
    Thus,
    \begin{equation*}
        \begin{cases}
            \crossentropyloss(P(y_2)y_1, P(y_2)y_2) & = H(\textnorm{y_2}^{-2}\textdotprod{y_1}{y_2}), \\
            \crossentropyloss(0, P(y_2)y_2) & = H(0),
        \end{cases}
        \quad \text{where}~
        \fct{H}
            {t}
            {\Rbb}
            {\mathrm{le}\left(ty_2\right) - t\textdotprod{y_2}{\softmax(y_2)}.}
            {\Rbb}
    \end{equation*}
    Moreover,~$H$ is twice differentiable, and for all~$t \in \Rbb$, we have~$H'(t) = \textdotprod{y_2}{\softmax(ty_2)-\softmax(y_2)}$ and~$H''(t) = \textdotprod{y_2^2}{\softmax(ty_2)}-\textdotprod{y_2}{\softmax(ty_2)}^2 = \mathrm{Var}_{i \sim \softmax(ty_2)}([y_2]_i) \geq 0$.
    So,~$H$ is constant when~$y_2$ is a multiple of~$\1$ and strictly convex with minimizer~$1$ otherwise.
    Then, for all~$t \in \Rbb$, we deduce that~$H(t) < H(0)$ implies~$t > 0$, and thus
    \begin{equation*}
        \crossentropyloss(P(y_2)y_1, P(y_2)y_2) < \crossentropyloss(0, P(y_2)y_2)
        \iff H(\textnorm{y_2}^{-2}\textdotprod{y_1}{y_2}) < H(0)
        \implies \textdotprod{y_1}{y_2} > 0.
    \end{equation*}
\end{remark}

\clearpage\newpage
\section{Hybrid Algorithm Leveraging \nn Attacks and \dfo}
\label{section:optimization_hybrid}

This section addresses the second goal of the paper, of designing a hybrid optimization algorithm that combines optimization through directional \nn attacks (as in Section~\ref{section:optimization_attack/use_in_optim}) with techniques from \emph{derivative-free optimization} (\dfo)~\cite{AuHa2017,CoScVibook,LaMeWi2019} to tackle Problem~\eqref{problem:P} regardless of the structure of~$\Phi$.
In this work, we focus to the \emph{Covering Direct Search Method (\cdsm)}~\cite{AuBoBo24Covering}.
This choice is motivated by its theoretical and practical reliability for problems with an objective function that is only continuous.

Section~\ref{section:optimization_hybrid/cdsm} reviews the \cdsm and its convergence to a local solution.
Then, Section~\ref{section:optimization_hybrid/hybrid} presents our hybrid Algorithm~\ref{algorithm:hybrid}, along with its convergence analysis (Theorem~\ref{theorem:convergence_hybrid}) inherited from the \cdsm.

Before going through this section, let us stress that we design our hybrid method from \dfo methods only to address our lack of assumptions about the structure of~$\Phi$.
The core idea of our methodology is to hybridize \nn attacks with an optimization algorithm that is well-suited for the problem at hand.
When~$\Phi$ is a glass-box \nn, it is preferable to consider an algorithm exploiting the explicit architecture of~$\Phi$, instead of \dfo methods oblivious to this structure.
Similarly, when~$\Phi$ is not a glass-box but the function it encodes is smooth, we expect model-based methods~\cite[Part~4]{AuHa2017} and Bayesian methods~\cite{GrRaGuVeVe20BayesianOptimization} to be competitive with the \cdsm.
We may even adapt the \bfgs algorithm~\cite{NoWr2006}, by replacing evaluations of the gradient of~$f \circ \Phi$ by the \attackstep operator.
Algorithm~\ref{algorithm:hybrid} may therefore consist in a baseline to assess the performance of more sophisticated methods.

\subsection{Direct Search Method and Covering Direct Search Method}
\label{section:optimization_hybrid/cdsm}

The \dsm proceeds iteratively, with each iteration comprising a \searchstep and a \pollstep steps.
The \searchstep step is optional, but it allows for user-defined strategies and globalization techniques.
The \pollstep step is mandatory and underpins all theoretical guarantees, at the cost of a more stringent definition.
The \dsm has been historically split into two subclasses: \emph{mesh \dsm}~\cite{AuHa2017} and \emph{sufficient increase \dsm}~\cite{CoScVibook}, each defining some specific restrictions on the \searchstep and \pollstep steps.
Both subclasses of \dsm ensure that some limit points they generate satisfy necessary optimality conditions for Problem~\eqref{problem:P}.
Nevertheless, an advance on \dsm unifies these two subclasses.
The \cdsm (\coveringstep \dsm)~\cite{AuBoBo24Covering} introduces a third step named the \coveringstep step, which allows for a stronger convergence analysis than \dsm.
The \coveringstep step, under mild assumptions, ensures convergence to local solutions (instead of only points satisfying necessary optimality conditions).
To our knowledge, at the time of writing the \cdsm is the only local \dfo algorithm with this convergence guarantee (although~\cite{AuBoBo24Covering} highlights that most \dfo algorithms may be enriched with a \coveringstep to obtain the same convergence properties).
We now overview the \searchstep and \pollstep steps as allowed in the \cdsm, and then we introduce the \coveringstep step.

The \searchstep step offers a permissive framework for evaluating many trial points.
It usually encloses globalization strategies such as the \emph{variable neighbourhood search}~\cite{HaMl01a,MlDrKoVuCa2008}.
At each iteration~$k \in \Nbb$, the \searchstep step usually relies on the current incumbent solution~$x^k$ and the \emph{trial points history}~$\Hcal^k$ (the set of all points evaluated by the algorithm up to iteration~$k$); and the set~$\searchstep(x^k,\Hcal^k) \subseteq \Rbb^n$ is only required to be finite.
Accordingly, we formalize a \emph{\searchstep operator} as follows.

\begin{definition}[\searchstep operator]
    \label{definition:search_operator}
    A set-valued function~$\searchstep: \Rbb^n \times 2^{\Rbb^n} \to 2^{\Rbb^n}$ is a \emph{\searchstep operator} if~$\searchstep(x,\Hcal) \subseteq \Rbb^n$ is empty or finite for all~$(x,\Hcal) \in \Rbb^n \times 2^{\Rbb^n}$.
\end{definition}

The \pollstep step has a more stringent definition.
At all iterations~$k \in \Nbb$, it evaluates points in the ball centred at~$x^k$ with radius~$r^k > 0$.
The literature usually defines the set~$\pollstep(x^k,r^k,\Hcal^k) \in \Bcal^n_{r^k}(x^k)$ of trial points so that the set~$\{t-x^k : t \in \pollstep(x^k,r^k,\Hcal^k)\} \in \Bcal^n_{r^k}$ of trial directions forms a positive spanning set of~$\Rbb^n$.
We formalize the \emph{\pollstep operator} accordingly.
A popular strategy~\cite{AbAuDeLe09} samples a unit vector~$v^k$, builds the matrix~$M^k \defequal I-2(v^k)\transpose{(v^k)}$, and sets~$\pollstep(x^k,r^k,\Hcal^k) \defequal \{x^k \pm r^kM^ke_i\}_{i=1}^{n}$.
In addition, popular approaches such as~\cite{CoLed2011,VeKa2017} use the set~$\Hcal^k$ to construct a local quadratic model of~$f \circ \Phi$ near~$x^k$ and add the gradient of that model at~$x^k$ to the set of trial directions.

\begin{definition}[\pollstep operator]
    \label{definition:poll_operator}
    A set-valued function~$\pollstep: \Rbb^n \times \Rbb_+ \times 2^{\Rbb^n} \to 2^{\Rbb^n}$ is a \emph{\pollstep operator} if for all~$(x,r,\Hcal) \in \Rbb^n \times \Rbb_+ \times 2^{\Rbb^n}$, the set~$\{t-x^k : t \in \pollstep(x,r,\Hcal)\} \subseteq \Bcal^n_r$ positively spans~$\Rbb^n$.
\end{definition}

Finally, the \coveringstep step~\cite{AuBoBo24Covering} may be added to the \dsm on top of the \searchstep and \pollstep steps.
This step consists, at all iterations~$k \in \Nbb$, in evaluating points in the ball~$\Bcal^n_1(x^k)$ that are far enough from the trial points history~$\Hcal^k$.
The fixed radius (here set to~$1$ for simplicity) ensures that all accumulation points of~$(x^k)_{k \in \Nbb}$ are local solutions.
Below, we consider a concise definition of the \coveringstep operator for simplicity, although we remark that~\cite{AuBoBo24Covering} formalizes others that allow for an easier computation.

\begin{definition}[\coveringstep operator]
    \label{definition:covering_operator}
    A set-valued function~$\coveringstep: \Rbb^n \times 2^{\Rbb^n} \to 2^{\Rbb^n}$ is a \emph{\coveringstep operator} if~$\coveringstep(x,\Hcal) \defequal \argmax_{x' \in \Bcal^n_1(x)} \dist(x',\Hcal)$ holds for all~$(x,\Hcal) \in \Rbb^n \times 2^{\Rbb^n}$ with~$\Hcal \neq \emptyset$.
\end{definition}

The \cdsm is summarized in Algorithm~\ref{algorithm:cdsm} below.
As we formalize in Proposition~\ref{proposition:convergence_cdsm}, its convergence properties established in~\cite[Theorem~1]{AuBoBo24Covering} hold for Problem~\eqref{problem:P}.

\begin{algorithm}[!htbp]
\caption{\cdsm to solve Problem~\eqref{problem:P}.}
\label{algorithm:cdsm}
\begin{algorithmic}

\IState{0} \algorithmicfunction \texttt{iteration}\_\cdsm(\coveringstep, \searchstep, \pollstep; $x$, $r$, $\Hcal$, $\tau$):
    \IState{1} \textbf{\coveringstep step}:
        \IState{2}
            set~$\TcalC \defequal \coveringstep(x, \Hcal)$;
        \IState{2} if~$\TcalC \cap F \neq \emptyset$,
            then select~$\tC \in \argmax f(\Phi(\TcalC \cap F))$,
            else set~$\tC \defequal x$;
        \IState{2} if~$f(\Phi(\tC)) > f(\Phi(x))$,
            then \algorithmicreturn $(\tC,~ r,~ \Hcal \cup \TcalC)$,
            else continue;
    \IState{1} \textbf{\searchstep step}:
        \IState{2}
            set~$\TcalS \defequal \searchstep(x, \Hcal)$;
        \IState{2} if~$\TcalS \cap F \neq \emptyset$,
            then select~$\tS \in \argmax f(\Phi(\TcalS \cap F))$,
            else set~$\tS \defequal x$;
        \IState{2} if~$f(\Phi(\tS)) > f(\Phi(x))$,
            then \algorithmicreturn $(\tS,~ r,~ \Hcal \cup \TcalC \cup \TcalS)$,
            else continue;
    \IState{1} \textbf{\pollstep step}:
        \IState{2}
            set~$\TcalP \defequal \pollstep(x, r, \Hcal)$;
        \IState{2} if~$\TcalP \cap F \neq \emptyset$,
            then select~$\tP \in \argmax f(\Phi(\TcalP \cap F))$,
            else set~$\tP \defequal x$;
        \IState{2} if~$f(\Phi(\tP)) > f(\Phi(x))$,
            then \algorithmicreturn $(\tP,~ \tau r,~ \Hcal \cup \TcalC \cup \TcalS \cup \TcalP)$,
            else continue;
    \IState{1} \textbf{\failurestep step}:
        \IState{2} \algorithmicreturn $(x,~ \tau^{-1}r,~ \Hcal \cup \TcalC \cup \TcalS \cup \TcalP)$.

\vspace{1em}
\IState{0} \textbf{Main algorithm}:
    \IState{1}
        select a~$\coveringstep$ operator,
        a~$\searchstep$ operator
        and a~$\pollstep$ operator;
        select~$\tau > 1$;
    \IState{1}
        select~$x^0 \in F$; select~$r^0 \in \Rbb_+^*$;
        set~$\Hcal^0 \defequal \{x^0\}$;
    \IState{1} \algorithmicfor~$k \in \Nbb$,
        set~$(x^{k+1}, r^{k+1}, \Hcal^{k+1}) \defequal \texttt{iteration}\_\cdsm(\coveringstep, \searchstep, \pollstep; x^k, r^k, \Hcal^k, \tau)$.

\end{algorithmic}
\end{algorithm}

\begin{proposition}[Adapted from~{\cite[Theorem~1]{AuBoBo24Covering}}]
\label{proposition:convergence_cdsm}
    Let~$(x^k)_{k \in \Nbb}$ be the sequence of incumbent solutions generated by Algorithm~\ref{algorithm:cdsm} solving Problem~\eqref{problem:P} under Assumption~\ref{assumption:problem}.
    Then~$(x^k)_{k \in \Nbb}$ admits at least one accumulation point, and all of them are local solutions to Problem~\eqref{problem:P}.
\end{proposition}

\subsection{A Hybrid Algorithm Using Directional \nn Attacks and \cdsm}
\label{section:optimization_hybrid/hybrid}

Algorithm~\ref{algorithm:cdsm} ignores the structure of Problem~\eqref{problem:P}.
Thus, it remains applicable to hard instances, but it usually converges slowly in practice.
In contrast, our \attackstep operator from Section~\ref{section:optimization_attack/use_in_optim}, relying on directional \nn attacks, offers a heuristic yet potentially efficient strategy for local improvements of the incumbent solution.
Then, each of these two approaches offsets the limitations of the other.
Motivated by this synergy, we propose a \emph{hybrid} algorithm that combines directional \nn attacks with the \cdsm.
At each iteration~$k \in \Nbb$, our algorithm first attempts a directional \nn attack via the \attackstep operator.
This attack provides a trial point~$x^k_\sharp$ that is tested against the incumbent solution~$x^k$.
Then the algorithm runs the \cdsm steps from the best of these two points.
This logic is formalized in Algorithm~\ref{algorithm:hybrid}.
We discuss some possible improvements of this algorithm in Remark~\ref{remark:hybrid_improvements}.

\begin{algorithm}[!ht]
\caption{Hybrid (directional \nn attacks with \cdsm) algorithm to solve Problem~\eqref{problem:P}.}
\label{algorithm:hybrid}
\begin{algorithmic}

\IState{0} \algorithmicfunction \texttt{iteration}\_\attackstep(\attackstep; $x$, $r$, $\Hcal$, $\tau$):
    \IState{1} \textbf{\attackstep step}:
        \IState{2}
            set~$\DcalA \defequal \attackstep(x, r)$;
            set~$\TcalA \defequal \left\{x + d,~ d \in \DcalA\right\}$;
        \IState{2} if~$\TcalA \cap F \neq \emptyset$,
            then select~$\tA \in \argmax f(\Phi(\TcalA \cap F))$,
            else set~$\tA \defequal x$;
        \IState{2} if~$f(\Phi(\tA)) > f(\Phi(x))$,
            then \algorithmicreturn $(\tA,~ \tau r,~ \Hcal \cup \TcalA)$,
            else continue;
    \IState{1} \textbf{\failurestep step}:
        \IState{2} \algorithmicreturn $(x,~ \tau^{-1}r,~ \Hcal \cup \TcalA)$.

\vspace{1em}
\IState{0} \textbf{Main algorithm}:

    \IState{1}
        select an~$\attackstep$ operator,
        a~$\coveringstep$ operator,
        a~$\searchstep$ operator
        and a~$\pollstep$ operator;

    \IState{1}
        select~$\tau_\atck > 1$ and~$\tau_\cdsm > 1$;

    \IState{1}
        select~$x^0 \in F$; select~$r_\atck^0 \in \Rbb_+^*$ and~$r_\cdsm^0 \in \Rbb_+^*$; set~$\Hcal^0 \defequal \{x^0\}$;

    \IState{1} \algorithmicfor~$k \in \Nbb$,
        \IState{2} set~$(x_\sharp^k, r_\atck^{k+1}, \Hcal_\sharp^k) \defequal \texttt{iteration}\_\attackstep(\attackstep; x^k, r_\atck^k, \Hcal^k, \tau_\atck)$;
        \IState{2} set~$(x^{k+1}, r_\cdsm^{k+1}, \Hcal^{k+1}) \defequal \texttt{iteration}\_\cdsm(\coveringstep, \searchstep, \pollstep; x_\sharp^k, r_\cdsm^k, \Hcal_\sharp^k, \tau_\cdsm)$.

\end{algorithmic}
\end{algorithm}

This algorithm inherits the convergence analysis of the \cdsm, as formalized in the next Theorem~\ref{theorem:convergence_hybrid}.

\begin{theorem}
    \label{theorem:convergence_hybrid}
    Let~$(x^k_\sharp)_{k \in \Nbb}$ be the sequence of incumbent solutions generated by Algorithm~\ref{algorithm:hybrid} solving Problem~\eqref{problem:P} under Assumption~\ref{assumption:problem}.
    Then~$(x^k_\sharp)_{k \in \Nbb}$ admits at least one accumulation point, and all such points are local solutions to Problem~\eqref{problem:P}.
\end{theorem}
\begin{proof}
    This proof relies on~\cite[Theorem~2]{AuBoBo24Covering}, which claims the following.
    Consider that Assumption~\ref{assumption:problem} holds, and define~$x^0 \in F$ and~$\Hcal^0 \subseteq \Rbb^n$.
    For all~$k \in \Nbb$, set~$\TcalC^k \defequal \coveringstep(x^k,\Hcal^k)$ and select~$\Hcal^{k+1}$ such that~$\Hcal^k \cup \TcalC^k \subseteq \Hcal^{k+1}$ and select~$x^{k+1} \in \argmax f(\Phi(\Hcal^{k+1} \cap F))$.
    Then~$(x^k)_{k \in \Nbb}$ admits at least one accumulation point, and all of them are local solutions to Problem~\eqref{problem:P}.

    Let~$(x^k_\sharp, \Hcal^k_\sharp)_{k \in \Nbb}$ be the sequence generated by Algorithm~\ref{algorithm:hybrid} under Assumption~\ref{assumption:problem}.
    By construction, for all~$k \in \Nbb$, it holds that~$\Hcal^k_\sharp \cup \coveringstep(x^k_\sharp, \Hcal^k_\sharp) \subseteq \Hcal^{k+1}_\sharp$, since~$\Hcal^{k+1}_\sharp$ equals~$\Hcal^k_\sharp$, plus all the trial points generated by the \cdsm steps from~$(x^k_\sharp, \Hcal^k_\sharp)$ (which include the \coveringstep step) and the \attackstep at iteration~$k+1$.
    Moreover, it also holds by design that~$x^{k+1}_\sharp \in \argmax f(\Phi(\Hcal^{k+1}_\sharp \cap F))$ for all~$k \in \Nbb$.
    Thus,~\cite[Theorem~2]{AuBoBo24Covering} is applicable to the sequence~$(x^k_\sharp, \Hcal^k_\sharp)_{k \in \Nbb}$, and our claim follows directly.
\end{proof}

To conclude this section, let us discuss the compatibility of Algorithm~\ref{algorithm:hybrid} with implementations of the \attackstep operator.
As we stress in Section~\ref{section:optimization_attack/use_in_optim/heuristic_tool}, current solvers for \nn attacks are incompatible with constrained instances of Problem~\eqref{problem:P}.
Our workaround may then yield attack directions that break the constraints.
However, Algorithm~\ref{algorithm:hybrid} is resilient to this case.
Indeed, all infeasible trial points generated by the \attackstep step are rejected at the definition of~$\tA$, so for all~$k \in \Nbb$, the incumbent solution~$x^k_\sharp$ after the \attackstep step remains feasible.
The worst case scenario (when the workaround always returns infeasible directions) is that the \attackstep never contributes to the optimization process.
Consequently, our workaround deteriorates neither the theoretical properties of Algorithm~\ref{algorithm:hybrid} nor the numerical reliability of the \cdsm steps.
More generally, the \attackstep step does not contribute to the theoretical properties of Algorithm~\ref{algorithm:hybrid}, even when the \attackstep operator satisfies Assumption~\ref{assumption:attack_guarantees_ascent_directions}, since Theorem~\ref{theorem:convergence_hybrid} does not require Assumption~\ref{assumption:attack_guarantees_ascent_directions}.

\begin{remark}
    \label{remark:hybrid_improvements}
    For ease of presentation, we write Algorithm~\ref{algorithm:hybrid} in a concise form.
    Nevertheless, Algorithm~\ref{algorithm:hybrid} is actually compatible as is with many refinements from \dfo, so several of its restrictions may be alleviated with no alterations in its theoretical properties.
    For example, we currently impose that~$x^0$ is feasible, but the \textit{progressive barrier}~\cite[Chapter~10]{AuHa2017} relax this requirement.
    Similarly, each iteration is expensive since it requires evaluating all the trial points in all steps, such as the~$n+1$ points (at least) in the \pollstep step given by a randomly chosen positive basis.
    It is actually possible to interrupt the evaluation of trial points in any of the steps at the first successful trial point.
    We may also use the history at all iterations to construct models of~$f \circ \Phi$, or trial directions based on momentum.
    Finally, it is also possible to skip entirely the \cdsm steps at some iterations~$k \in \Nbb$, when the \attackstep step identifies an incumbent solution~$x^k_\sharp$ that passes a \emph{sufficient increase} check over~$x^k$, determined by~$f(\Phi(x^k_\sharp)) \geq f(\Phi(x^k)) + \tau\textabs{f(\Phi(x))}+\varepsilon$ with~$(\tau,\varepsilon) \in \Rbb_+ \times \Rbb_+^*$.
\end{remark}

\section{Numerical Experiments}
\label{section:numerical}

In this section, we evaluate the numerical performance and behaviour of Algorithm~\ref{algorithm:hybrid} on three problems with diverse structures.
On each problem, we conduct the next three analyses.

\paragraph{Experiment 1.}
\label{experiment:performance_algo}
\textbf{General performance comparison.}
We evaluate the performance of Algorithm~\ref{algorithm:hybrid} against four baselines: the \emph{steepest ascent algorithm} and the \emph{Broyden-Fletcher-Goldfarb-Shanno} (\bfgs) algorithm, both estimating gradients by backpropagations, the \cdsm, and the \attackstep operator alone.
For each algorithm, we track the number of calls to the \nn~$\Phi$ (either in forward and backward passes) and the elapsed runtime, and we display the best objective value found as functions of these metrics.
We repeat so five times with five different seeds for the random number generator.

\paragraph{Experiment 2.}
\label{experiment:contribution_steps}
\textbf{Contribution of each step.}
We analyze the respective roles of the \attackstep and \cdsm steps within Algorithm~\ref{algorithm:hybrid}.
At each iteration, we report two flags that describe which step of, respectively, the \texttt{iteration}\_\attackstep function and the \texttt{iteration}\_\cdsm function triggers its \algorithmicreturn statement.
For comparison, we also conduct the same analysis for Algorithm~\ref{algorithm:cdsm} using only the \cdsm.

\paragraph{Experiment 3.}
\label{experiment:potential_attack}
\textbf{Alternative \attackstep operators.}
We evaluate the effectiveness of several~$\attackstep_{(\Srm,\lambda)}$ operators as proposed in Section~\ref{section:optimization_attack/use_in_optim/heuristic_tool}.
We define three such operators from the Torchattacks~\cite{kim2020torchattacks} versions of the \ffgsm~\cite{WoRiKo20FFGSM} and \rfgsm~\cite{TrKuPaGoBoMcD20RFGSM} and \pgd~\cite{Madry19PGD} algorithms respectively, and one operator based on backpropagation (returning~$d \defequal b \ r/\textnorm{b}$ with~$b \defequal \backprop(\Phi, x^{k(j)}, \nabla f(\Phi(x^{k(j)})))$).
To this end, we extract a sample~$(x^{k(j)})_{j \in \Jbb}$ (with~$\Jbb \subset \Nbb$) of incumbent solutions from the sequence~$(x^k)_{k \in \Nbb}$ generated by the \cdsm.
This sample spans a range of objective values for Problem~\eqref{problem:P}, enabling evaluation across different optimization stages.
For each~$j \in \Jbb$ and several radii~$r \in \Rbb_+$, we test whether~$\attackstep_{(\Srm,\lambda)}(x^{k(j)}, r)$ yields an ascent direction~$d$.
We then mark by a bullet the case where~$d$ yields an increase of at least~$1\%$ of the optimality gap between~$f(\Phi(x^{k(j)}))$ and the solution returned by the \cdsm, and by a cross the case where~$d$ yields an increase smaller than this threshold.
We also track the runtime and the associated values~$f(\Phi(x^{k(j)}))$ to relate ascent potential with point quality.

The full numerical setup is presented in Section~\ref{section:numerical/setting}.
Section~\ref{section:numerical/barycenter_into_resnet} explores a task to align the output of the \resnet network (introduced in Section~\ref{section:optimization_attack/notion/origin_illustration}) with a fixed target.
This is a proof-of-concept task, which offers a trade-off between a \nn naturally sensitive to attacks and a challenging optimization landscape.
Then, Section~\ref{section:numerical/warcraft} addresses a counterfactual explanation task involving a \nn that generates Warcraft maps, as proposed in~\cite{VAFoPaVi24CFOPT}.
This problem shows the performance of our method in a problem where the \nn involves a \emph{variational auto-encoder}, an architecture useful for optimization since it aggregates high-dimensional data into a low-dimensional manifold.
Finally, Section~\ref{section:numerical/bio_pinn} tackles a chemical engineering optimization problem of maximizing the production of some bio-diesel, using a \emph{Physics-Informed \nn} (\pinn) from~\cite{BiBoBl24PINNbiodiesel} which models a chemical reaction based on reaction time and the power of a heat flux.
This last problem highlights a context where an over-reliance on gradient-based ascent directions embarrasses the solution process.

Overall, our experiments show that Algorithm~\ref{algorithm:hybrid} is competitive with the baselines.
It outperforms them in the first two problems.
On the third problem, Algorithm~\ref{algorithm:hybrid} is dominated by the \cdsm baseline in terms of elapsed time but dominates in terms of number of passes through~$\Phi$.
It also appears that the \cdsm steps consistently repositions the incumbent solution when the \attackstep struggles, so that the \attackstep step contributes to the performance of Algorithm~\ref{algorithm:hybrid} across all optimization stages.
Finally, all~$\attackstep_{(\Srm,\lambda)}$ operators based on attack algorithms have similar potential.
However, the \ffgsm and \rfgsm algorithms are sensibly faster than the \pgd algorithm.
We also observe that the~$\attackstep_{(\Srm,\lambda)}$ operator based on attack algorithms outperform the backpropagation to identify ascent direction.

\subsection{Experimental Setup}
\label{section:numerical/setting}

As discussed in Sections~\ref{section:optimization_attack/use_in_optim/heuristic_tool} and~\ref{section:optimization_hybrid/hybrid}, current solvers for \nn attacks are incompatible with our setting studied in Section~\ref{section:optimization_attack/use_in_optim/setting_ascent_directions}.
This does not impact the properties of our hybrid Algorithm~\ref{algorithm:hybrid}.
Yet, one of our baselines consists in repeated calls to the $\attackstep_{(\Srm,\lambda)}$ operator from Definition~\ref{definition:attack_operator_practical}, which solves an unconstrained relaxed attack problem.
Consequently, for consistency in our experiments, we run all methods on the relaxation of Problem~\eqref{problem:P} obtained similarly,
\begin{equation*}
    \label{problem:P_practical}
    \tag{$\widetilde{\Pbf}$}
    \problemoptimfree{\maximize}{x \in \Rbb^n}{\widetilde{f}(\widetilde{\Phi}(x)),}
\end{equation*}
which replaces the constraint~$c(\Phi(x)) \leq 0$ by a quadratic penalization (with weight~$\lambda \in \Rbb_+$) via
\begin{equation*}
    \fct
        {\widetilde{\Phi}}
        {x}
        {\Rbb^n}
        {\begin{bmatrix} \Phi(x) \\ \relu(c(\Phi(x))) \end{bmatrix}}
        {\Rbb^{m+p}}
    \qquad \mbox{and} \qquad
    \fct
        {\widetilde{f}}
        {\begin{bmatrix} y \\ z \end{bmatrix}}
        {\Rbb^{m+p}}
        {f(y)-\lambda\norm{z}_2^2.}
        {\Rbb}
\end{equation*}
The solutions to Problem~\eqref{problem:P_practical} may differ from those of Problem~\eqref{problem:P}.
For each problem, we set~$\lambda$ to a small value to limit the infeasibility of the solutions.
Nonetheless, for each method with returned solution~$x^*$, we measure the feasibility of~$x^*$ in Problem~\eqref{problem:P} by calculating~$c(\Phi(x^*))$.
The methods we compare in our experiments are described below.
We fine-tuned their parameters using numerical evidence gathered during preliminary tests and similarity with common practice from \dfo.

\newcommand{\methodatck}{\Mbb_{\atck}}
\paragraph{Method~$\methodatck$: \attackstep operator only.}
This method uses the~$\attackstep_{(\Srm,\lambda)}$ operator from Definition~\ref{definition:attack_operator_practical}, instantiated with the Torchattacks~\cite{kim2020torchattacks} implementation of the \rfgsm algorithm~\cite{TrKuPaGoBoMcD20RFGSM} under~$\textnorm{\cdot}_\infty$ and~$\crossentropyloss$ loss function.
We conducted preliminary experiments with a variant relying on the \pgd algorithm~\cite{Madry19PGD}, but we observed that this stronger attack yields similar results, so we favour \rfgsm since it is faster.
Given~$x^0 \in \Rbb^n$ and~$r^0 \in \Rbb_+^*$ and~$\Hcal^0 \defequal \{x^0\}$, and an expansion parameter~$\tau \defequal 1.1$, the algorithm iterates for all~$k \in \Nbb$ as~$(x^{k+1}, r^{k+1}, \Hcal^{k+1}) \defequal \texttt{iteration}\_\attackstep(\attackstep_{(\Srm,\lambda)}; x^k, r^k, \Hcal^k, \tau)$.

\newcommand{\methodbprp}{\Mbb_{\bprp}}
\paragraph{Method~$\methodbprp$: first-order method using backpropagation.}
This method works similarly to the~$\methodatck$ method, except that it re-defines the~$\attackstep_{(\Srm,\lambda)}$ operator so that, for all~$(x,r) \in \Rbb^n \times \Rbb_+$, the returned set of directions is the singleton~$\attackstep(x,r) \defequal \{rd\textnorm{d}^{-1}\}$ where~$d \defequal \backprop(\widetilde{\Phi}, x, \nabla \widetilde{f}(\widetilde{\Phi}(x)))$.

\newcommand{\methodbfgs}{\Mbb_{\bfgs}}
\paragraph{Method~$\methodbfgs$: second-order method  using backpropagation.}
This method applies the PyTorch~\cite{PyTorch} implementation of the \bfgs algorithm, with~$x \mapsto \backprop(\widetilde{\Phi}, x, \nabla \widetilde{f}(\widetilde{\Phi}(x)))$ as a surrogate of the composite gradient function~$x \mapsto \nabla(\widetilde{f} \circ \widetilde{\Phi})(x)$.

\newcommand{\methodcdsm}{\Mbb_{\cdsm}}
\paragraph{Method~$\methodcdsm$: \cdsm baseline.}
This method runs Algorithm~\ref{algorithm:cdsm} on Problem~\eqref{problem:P_practical}.
For each~$k \in \Nbb$, the \coveringstep step is simplified so it evaluates a random point in the ball~$\Bcal_1(x^k)$ chosen uniformly (this does not match Definition~\ref{definition:covering_operator}, but~\cite{AuBoBo24Covering} proves that this preserves Proposition~\ref{proposition:convergence_cdsm} in practice), the \pollstep step follows the scheme from~\cite{AbAuDeLe09}, and the \searchstep step is skipped.
We choose~$\tau \defequal 2$, as usual in \dfo.

\newcommand{\methodhybr}{\Mbb_{\hybr}}
\paragraph{Method~$\methodhybr$: hybrid method.}
This method runs Algorithm~\ref{algorithm:hybrid} on Problem~\eqref{problem:P_practical}.
We define~$\attackstep_{(\Srm,\lambda)}$ from the Torchattacks~\cite{kim2020torchattacks} implementation of the \ffgsm~\cite{WoRiKo20FFGSM} algorithm, with~$\tau_\atck \defequal 1.2$, and the \cdsm steps and parameters as in the~$\methodcdsm$ method.

Each algorithm terminates when the runtime exceeds~$900$ seconds on the first problem and~$300$ seconds on the last two, or when all radii become smaller than~$10^{-5}$.
We implement our methods in Python~$3.12.2$, running on a single-thread Intel Xeon Gold 6258R CPU cadenced at~$2.70$GHz.
The code is available at \url{https://github.com/PierreYvesBouchet/Optimization-by-Directional-Attacks}.
In this section, we report only the main numerical evidence we gathered during our experiments.
We also report results for the~$\methodatck$ and~$\methodhybr$ methods using the same attack algorithm for the $\attackstep_{(\Srm,\lambda)}$ operator in all problems.
The algorithm we choose is a trade-off resulting from our results in a complete numerical campaign where, for both methods, we tested several \attackstep algorithms on each problem.
The best attack algorithm to construct the~$\attackstep_{(\Srm,\lambda)}$ operator is actually problem-dependent, and on each problem, we may obtain better results than those reported below by using a different \attackstep operator.
We provide the complete experimental materials in the GitHub repository.

\subsection{Proof-of-concept problem: Target Image Recovery}
\label{section:numerical/barycenter_into_resnet}

This task builds upon the \resnet classifier discussed in Section~\ref{section:optimization_attack/notion/origin_illustration} and on the \emph{barycentric image} function~$\bary: \Ibb^n \times \Rbb^n \to \Ibb$ defined by~$\bary(I;w) \defequal \sum_{\ell=1}^{n} w_\ell I_\ell$ for all sets~$I \defequal (I_\ell)_{\ell=1}^{n} \in \Ibb^n$ of~$n$ images and all vectors~$w \in \Pbb^n$ of weights.
We fix~$n \defequal 100$, and we select~$I \in \Ibb^n$ such that the input images are mutually dissimilar.
Then, we define
\begin{equation*}
    \fct
        {\Phi}
        {x}
        {\Rbb^n}
        {\resnet(\bary(I;\softmax(x)))}
        {\Pbb^{1000}}
    \quad \mbox{and} \quad
    \fct
        {f}
        {y}
        {\Pbb^{1000}}
        {-\norm{y-\resnet(I_1)}}
        {\Rbb}
\end{equation*}
and~$c: x \mapsto \begin{bmatrix} (x_i-10)(x_i+10) \end{bmatrix}_{i=1}^{n}$, defining the feasible set~$F \defequal [-10, 10]^n$.
The solution~$x^* \in \Rbb^n$ is the vector with~$x^*_1 \defequal 10$ and~$x^*_i \defequal -10$ for all~$i \neq 1$, and we initialize~$x^0 \defequal 0$.
We set~$\lambda \defequal 10$ for the penalization in the relaxed problem.
This problem is realistic in terms of dimension and complexity of the landscape of the objective function.
In particular, escaping the neighbourhood of the starting point is challenging for gradient-based methods.
We propose this problem for experiments on a controlled instance, since it involves a \nn that is known to be sensitive to attacks.

\begin{figure}[!ht]
    \centering
    \includegraphics[width=\linewidth]{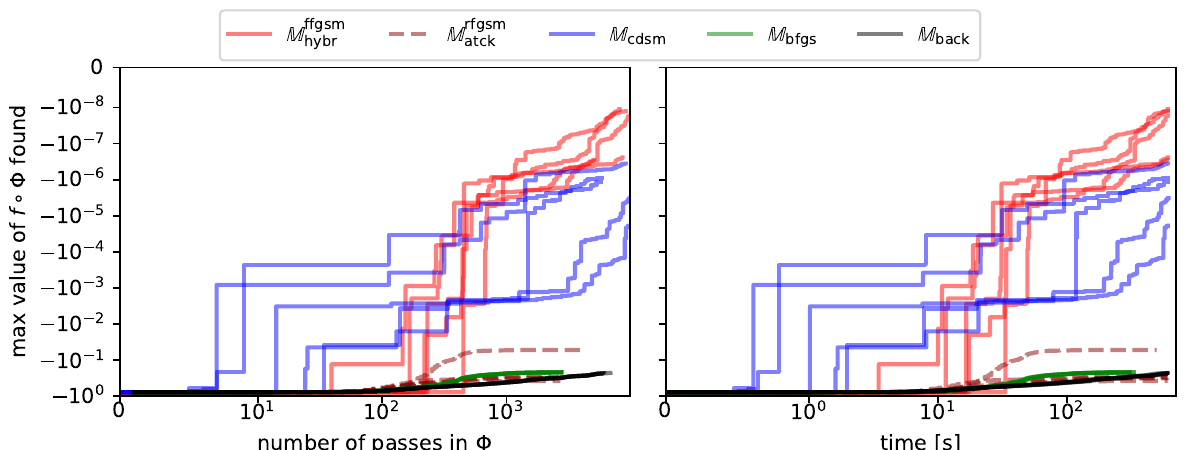}
    \caption{\nameref{experiment:performance_algo} in the Target Image Recovery problem from Section~\ref{section:numerical/barycenter_into_resnet}. The global solution to the problem has an objective value of~$-2.7E^{-13}$.}
    \label{figure:resnet/results_nb_pts}
\end{figure}

Let us analyze first the results of~\nameref{experiment:performance_algo} on this problem.
Figure~\ref{figure:resnet/results_nb_pts} shows that among the four baselines, only the methods involving \dfo converge to a near-optimal solution.
Indeed, the~$\methodatck$ method underperforms, as well as the first-order and second-order baselines~$\methodbprp$ and~$\methodbfgs$.
In contrast, the~$\methodcdsm$ method succeeds in approximating closely the optimum.
The~$\methodhybr$ method has a slower start than the~$\methodcdsm$ method due to the lack of performance of the \attackstep step near the starting point, but the \attackstep step eventually performs better when the algorithm escapes this challenging area.
As a result, the $\methodhybr$ method eventually dominates the~$\methodcdsm$ method by several orders of magnitude.
This experiment highlights how gradient-free method allows escaping the challenging area near the starting point.
All runs of all methods return feasible solutions.

\begin{figure}[!ht]
    \centering
    \includegraphics[width=\linewidth]{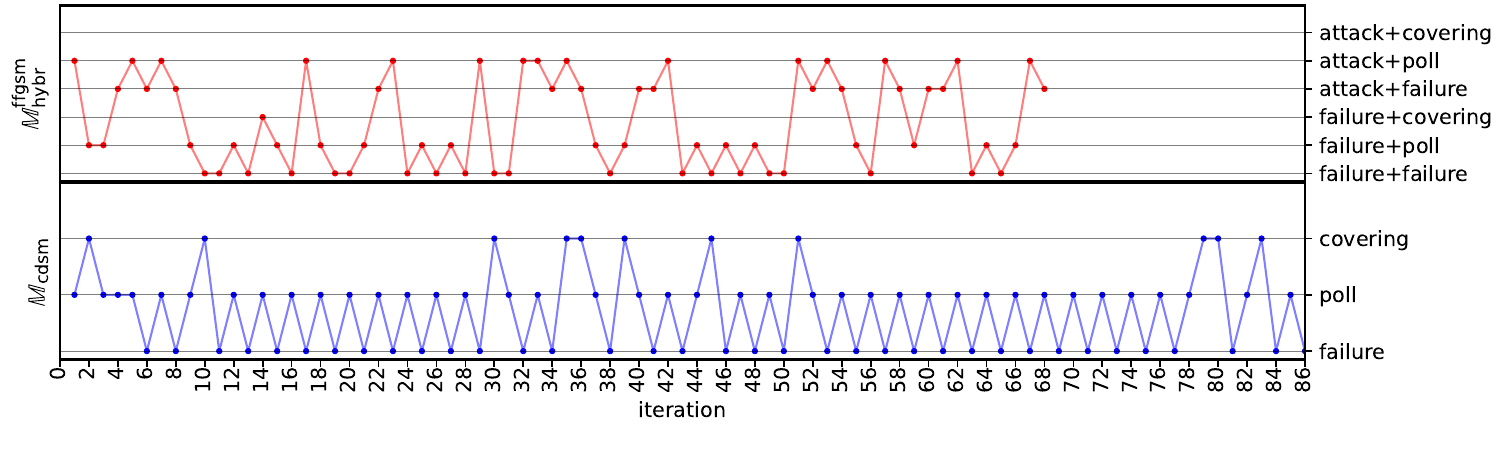}
    \caption{\nameref{experiment:contribution_steps} in the Target Image Recovery problem from Section~\ref{section:numerical/barycenter_into_resnet}.}
    \label{figure:resnet/iterations_results}
\end{figure}

We now focus to~\nameref{experiment:contribution_steps}, depicted in Figure~\ref{figure:resnet/iterations_results}.
The~$\methodcdsm$ method mostly progresses with the \pollstep step and a noticeable failure-success pattern.
With a detailed analysis to explain so, we observed that early successes increase the \pollstep radius until it reaches a value of~$16$, too large for the bounded feasible domain.
From that point, the repeated pattern is that one iteration fails, so the radius shrinks to~$8$, and then the following iteration successes, so the radius grows again to~$16$.
This domination of the \pollstep step among the \cdsm components carries over to the~$\methodhybr$ method.
However, noticeably, the performance of the~$\methodhybr$ method also results from the \attackstep step, which succeeds roughly half of the time.
This highlights the strength of combining attack techniques with a local search.
We believe that the \pollstep step assists the \attackstep step by repositioning the point at which the attack is more likely to perform.
This explains how the~$\methodhybr$ method escapes the neighbourhood of~$x^0$, which appears challenging for the \attackstep step.

\begin{figure}[!ht]
    \centering
    \includegraphics[width=\linewidth]{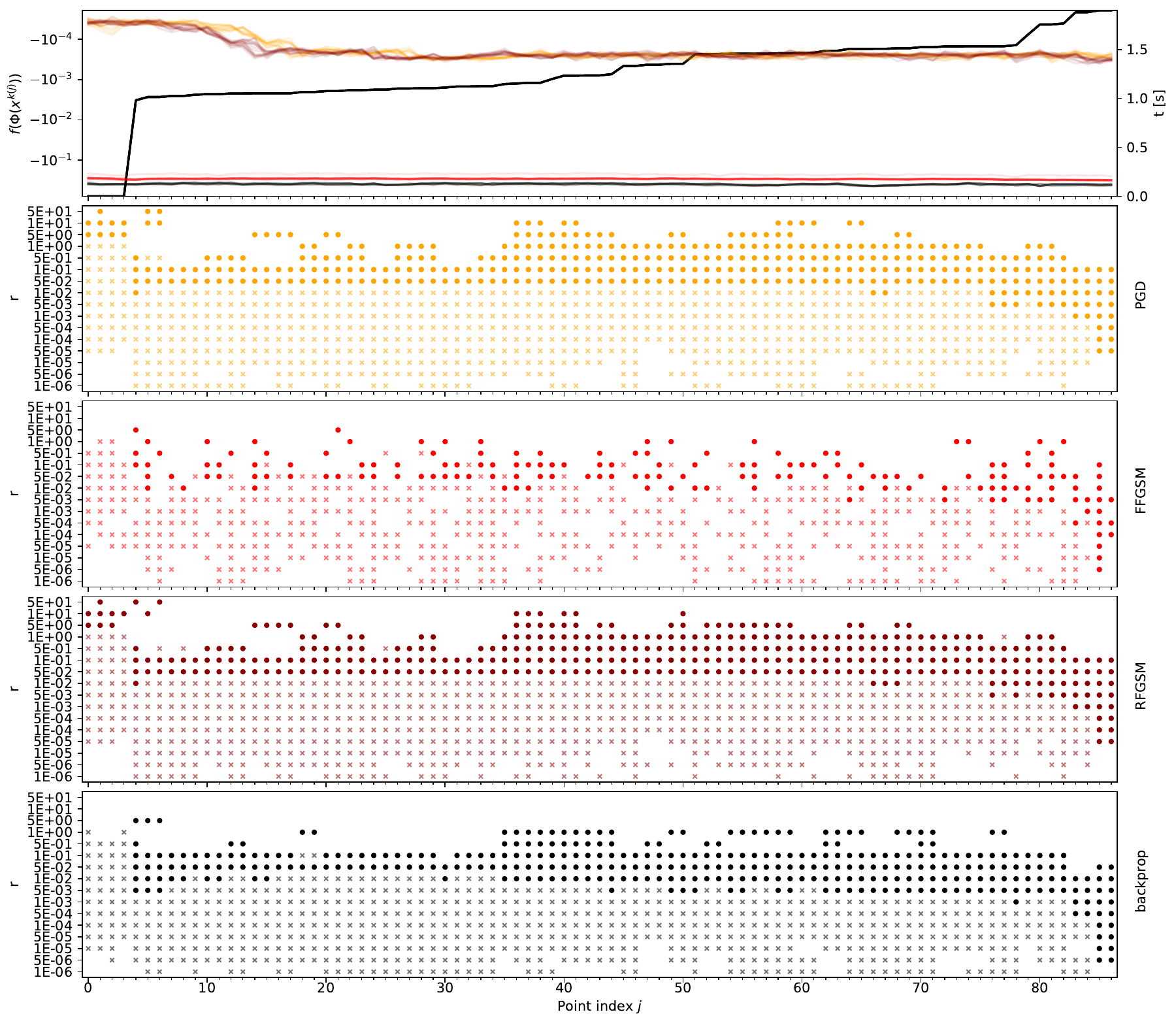}
    \caption{\nameref{experiment:potential_attack} in the Target Image Recovery problem from Section~\ref{section:numerical/barycenter_into_resnet}.}
    \label{figure:resnet/attacks_analysis}
\end{figure}

Finally,~\nameref{experiment:potential_attack} confirms the reliability of the \attackstep operator in identifying ascent directions throughout the optimization process.
Figure~\ref{figure:resnet/attacks_analysis} shows that \rfgsm and \pgd are more effective than backpropagation, especially in the regime where~$r \geq 1E{-}01$.
We also observe that on this problem, \ffgsm is as fast as the backpropagation, although less efficient.

\subsection{Application: Counterfactual Warcraft Maps}
\label{section:numerical/warcraft}

Our second problem, adapted from~\cite{VAFoPaVi24CFOPT}, focuses on \emph{optimal counterfactual explanations for structured prediction models} involving a \nn providing parameters of a combinatorial optimization layer.
We adapt the flagship application discussed in~\cite[Appendix~C]{VAFoPaVi24CFOPT}, which searches for \emph{$\varepsilon$-relative counterfactual Warcraft maps} with respect to the lightest path problem, in the sense defined in Section~\ref{section:introduction}.

\begin{figure}[!ht]
    \centering
    \includegraphics[width=0.75\linewidth]{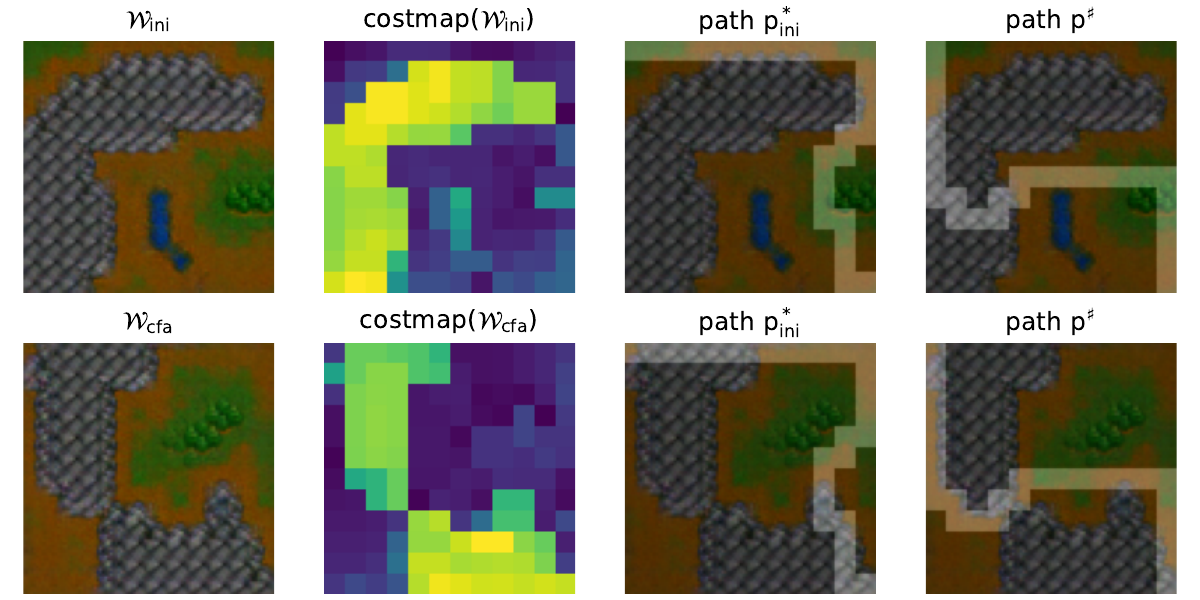}
    \caption{
        (First line) Warcraft map~$\Wcalini$, its associated~$\costmap$ output, the lightest path~$\pathmapini$ to reach the South-East from the North-West, and an alternative path~$\pathmap^\sharp$.
        (Second line) Similar displays for a counterfactual map~$\Wcal_\counterfactual$ with respect to~$\pathmap^\sharp$.
        This counterfactual is likely not optimal, since the two maps share limited similarities besides the surface of the mountainous area.
        Nevertheless,~$\Wcal_\counterfactual$ is well suited for~$\pathmap^\sharp$ since the mountain has a gorge exactly where~$\pathmap^\sharp$ crosses.
        }
    \label{figure:warcraft/example_maps}
\end{figure}

The video game \emph{Warcraft} is a tactical board game where some areas of the maps (represented as images in the space~$\Wbb \defequal [0,1]^{3 \times 96 \times 96}$ endowed with some norm~$\textnorm{\cdot}_\Wbb$) are lighter to cross than others.
To reflect this,~\cite{VAFoPaVi24CFOPT} designs a \nn~$\costmap: \Wbb \to \Cbb \defequal \Rbb_+^{12 \times 12}$ that splits any map~$\Wcal \in \Wbb$ into~$12 \times 12$ areas of~$8 \times 8$ pixels and estimates the cost for a character to enter these areas.
We endow the space~$\Cbb$ with the sum-of-entries norm~$\textnorm{\cdot}_\Cbb$ and the entry-wise inner product~$\odot$.
Thus, for any Warcraft map~$\Wcal \in \Wbb$ and any path~$\pathmap \in \{0,1\}^{12 \times 12}$ crossing the map, the cost of~$\pathmap$ along~$\Wcal$ reads as
\begin{equation*}
    \costpath(\Wcal;\pathmap) \defequal \norm{\costmap(\Wcal) \odot \pathmap}_\Cbb.
\end{equation*}
We consider a Warcraft map~$\Wcalini \in \Wbb$ and its cost map~$\Ccalini \defequal \costmap(\Wcalini) \in \Cbb$, and we compute the lightest path~$\pathmapini \in \{0,1\}^{12 \times 12}$ joining the Northwestern corner of~$\Wcalini$ to the Southeastern one.
Then we pick an alternative path~$\pathmap^\sharp \neq \pathmapini$ joining the same corners.
Our goal is to produce an alternative map that remains close to~$\Wcalini$ and favours this alternative path over the original one.
More specifically, we seek an \emph{$\varepsilon$-relative counterfactual explanation of~$\Wcalini$ with respect to the function~$\costpath$ and the path~$\pathmap^\sharp$}, where we fix~$\varepsilon \defequal \frac{1}{2}$.
We therefore search for a map~$\Wcal_\counterfactual \in \Wbb$ that solves the problem
\begin{equation*}
    \problemoptimoneline
        {\minimize}
        {\Wcal \in \Wbb}
        {\norm{\Wcal-\Wcalini}_\Wbb^2}
        {(1+\varepsilon)\costpath(\Wcal;\pathmap^\sharp) \leq \costpath(\Wcal;\pathmapini)}.
\end{equation*}
In other words, the goal is to find an alternative map~$\Wcal_\counterfactual \in \Wbb$ as close as possible to~$\Wcalini$, but on which the path~$\pathmap^\sharp$ is at least~$1.5$ times lighter than the path~$\pathmapini$.
We illustrate this problem in Figure~\ref{figure:warcraft/example_maps}.

However, the above problem has no constraints guaranteeing that the image encoded by~$\Wcal_\counterfactual$ is visually similar to~$\Wcalini$, or more generally to a real map from the game.
To enforce this requirement, we rely on another \nn from~\cite{VAFoPaVi24CFOPT}, that we denote by~$\warcraft: \Xbb \to \Wbb$ (where~$\Xbb \defequal \Rbb^n$ with~$n \defequal 64$).
This \nn is designed to generate credible Warcraft maps from abstract input vectors, where credible means that the images generated by~$\warcraft$ may not represent maps that truly exist in-game, but that look like real ones.
By design, all~$x \in \Rbb^n$ with~$\textabs{\textnorm{x}_\Xbb-\sqrt{n}} \leq 1$ yield~$\warcraft(x)$ being credible.
We consider that we know the vector~$\xini \in \Rbb^n$ that generates~$\Wcalini$ as~$\Wcalini \defequal \warcraft(\xini)$.
By design, any~$x \approx \xini$ generates~$\Wcal \defequal \warcraft(x) \approx \Wcalini$, so~$\textnorm{x-\xini}_\Xbb$ acts as a surrogate of~$\textnorm{\Wcal-\Wcalini}_\Wbb$.
Moreover, to enforce even further proximity between maps, we consider~$\textnorm{\costmap(\Wcal)-\Ccalini}_\Cbb$ as another surrogate of~$\textnorm{\Wcal-\Wcalini}_\Wbb$.
Then, instead of searching for~$\Wcal_\counterfactual \in \Wbb$ directly, we seek a counterfactual~$x_\counterfactual \in \Xbb$ to then obtain~$\Wcal_\counterfactual \defequal \warcraft(x_\counterfactual)$.
We also use the surrogate norms instead of the true norm.
We therefore solve
\begin{equation*}
    \problemoptim
        {\minimize}
        {x \in \Xbb}
        {\norm{x-\xini}_\Xbb^2 + \norm{\costmap(\warcraft(x))-\Ccalini}_\Cbb^2}
        {\norm{x}_\Xbb \in \left[\sqrt{n}-1 , \sqrt{n}+1\right], \\
         (1+\varepsilon)\costpath(\warcraft(x);\pathmap^\sharp) \leq \costpath(\warcraft(x);\pathmapini).
        }
\end{equation*}
This reformulation is easier to solve than the original problem, since~$\Xbb$ is a usual vector space instead of a space of images, and its dimension~$n = 64$ is much smaller than those of~$\Wbb$ ($3 \times 96 \times 96 = 27648$).
Moreover, for any~$x \in \Xbb$, we do not need to record~$\warcraft(x)$.
Either the objective and the constraint involving the~$\costpath$ function may be expressed in terms of~$\costmap(\warcraft(x))$ directly.

We fit this problem in the framework of Problem~\eqref{problem:P} as follows.
First, the \nn~$\Phi$ is
\begin{equation*}
    \fct
        {\Phi}
        {x}
        {\Xbb}
        {\left(x, \costmap\left(\warcraft\left(x\right)\right) \right)}
        {\Ybb \defequal \Xbb \times \Cbb}
\end{equation*}
and the goal function~$f$ is given by
\begin{equation*}
    \fct
        {f}
        {y \defequal (x,c)}
        {\Ybb}
        {-\left( \norm{x-\xini}_\Xbb^2 + \norm{c-\Ccalini}_\Cbb^2\right).}
        {\Rbb}
\end{equation*}

\begin{figure}[!ht]
    \centering
    \includegraphics[width=0.9\linewidth]{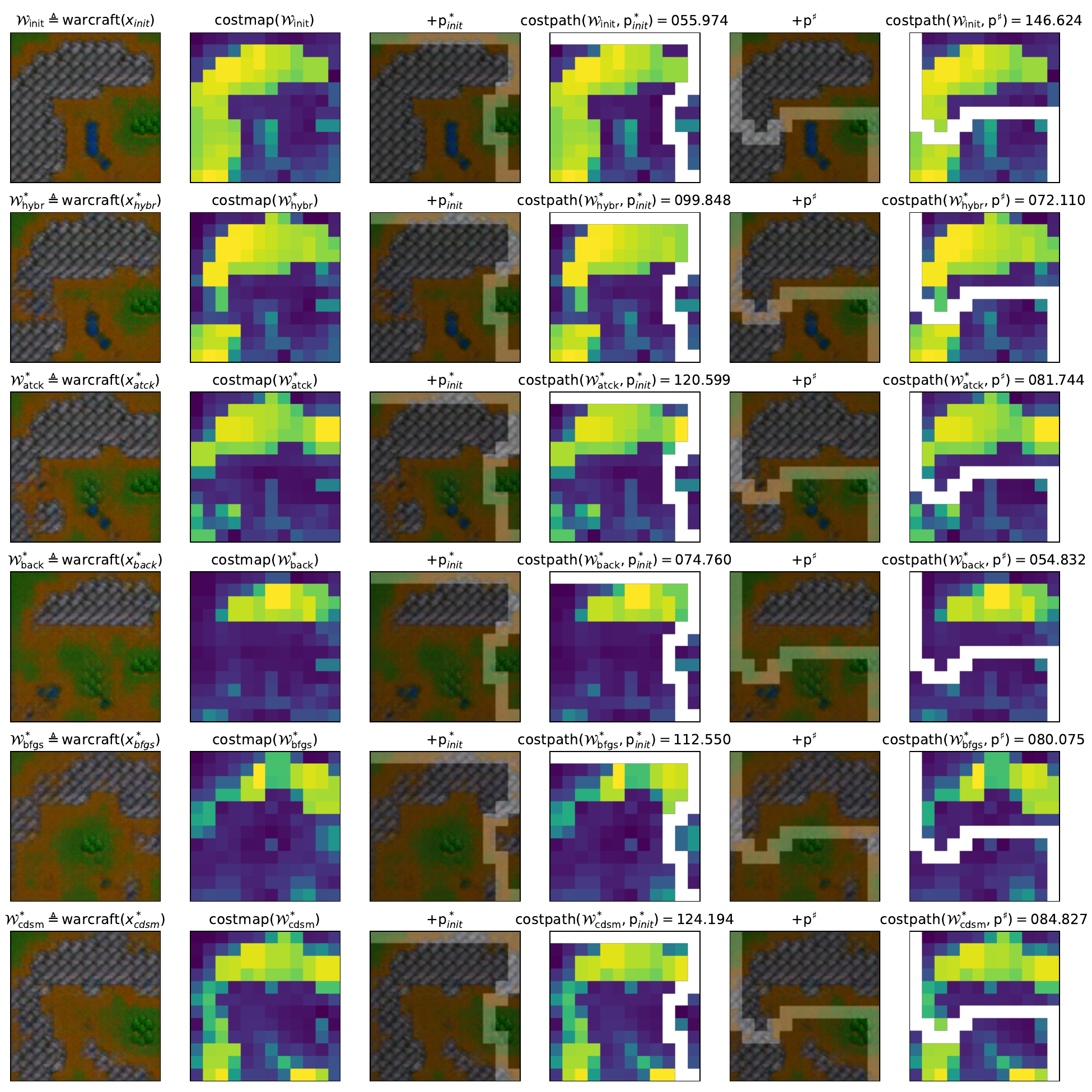}
    \caption{
        Warcraft maps considered in Section~\ref{section:numerical/warcraft}.
        The first line is related to~$x \defequal \xini$.
        Other lines are related to~$x$ being the solution returned by each method we test in our experiments.
        Each group of two consecutive columns displays~$\warcraft(x)$ and~$\costmap(\warcraft(x))$, then a visualization of the paths~$\pathmapini$ and~$\pathmap^\sharp$ and their costs.
        }
    \label{figure:warcraft/results_maps}
\end{figure}

Next, we define the constraints function by
\begin{equation*}
    \fct
        {c}
        {y \defequal (x,c)}
        {\Ybb}
        {\begin{bmatrix}
            \left(\norm{x}_\Xbb-\sqrt{n}-1\right)\left(\norm{x}_\Xbb-\sqrt{n}+1\right) \\[0.25em]
            (1+\varepsilon)\norm{c \odot \pathmap^\sharp}_\Cbb - \norm{c \odot \pathmapini}_\Cbb
            \end{bmatrix},
        }
        {\Rbb^2}
\end{equation*}
so that all~$x \in \Xbb$ with~$c(\Phi(x)) \leq 0$ yield a credible and valid counterfactual explanation as~$\warcraft(x)$; since~$\textabs{\textnorm{x}_\Xbb-\sqrt{n}} \leq 1$ and~$(1+\varepsilon)\costpath(\warcraft(x),\pathmap^\sharp) \leq \costpath(\warcraft(x),\pathmapini)$.
We initialize all methods with~$x^0 \defequal \xini$, and~$\lambda \defequal 10$.
The counterfactual maps calculated by each method are shown in Figure~\ref{figure:warcraft/results_maps}.
Interestingly, the map returned by the~$\methodhybr$ method is visually satisfying because of its proximity to the initial map; the only sensible modification being a tight gorge through the mountain to lighten the path~$\pathmap^\sharp$.
All methods return solutions that break only the constraint related to the~$\varepsilon$-relative counterfactual.
In the worst case, the returned map is an~$\varepsilon$-relative counterfactual with~$\varepsilon_\hybr \approx 1.42$,~$\varepsilon_\atck \approx 1.48$,~$\varepsilon_\bprp \approx 1.36$,~$\varepsilon_\bfgs \approx 1.41$, and~$\varepsilon_\cdsm \approx 1.46$.

\begin{figure}[!ht]
    \centering
    \includegraphics[width=0.9\linewidth]{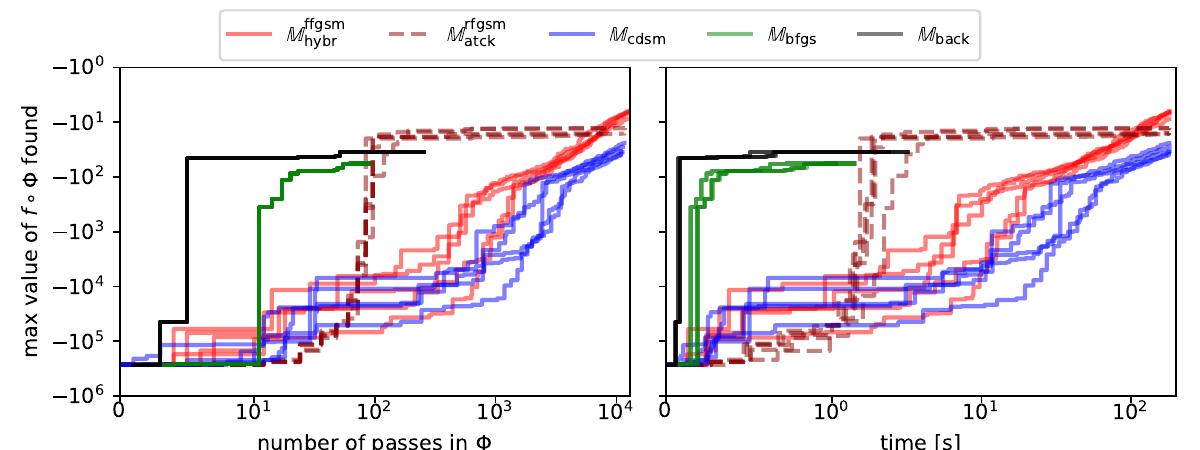}
    \caption{\nameref{experiment:performance_algo} in the Counterfactual Warcraft Maps problem from Section~\ref{section:numerical/warcraft}.}
    \label{figure:warcraft/results_nb_pts}
\end{figure}

\nameref{experiment:performance_algo}, displayed in Figure~\ref{figure:warcraft/results_nb_pts}, shows that the~$\methodhybr$ method delivers the best final output.
The~$\methodatck$ and~$\methodbprp$ and~$\methodbfgs$ methods quickly approach an almost feasible solution but fail to refine it further.
Although~$\methodbprp$ and~$\methodbfgs$ converge faster than~$\methodatck$ in this problem, their returned solutions violate the counterfactual constraint by a sensibly larger amount than the solution returned by the~$\methodatck$ method.
Our hybrid method~$\methodhybr$ acts as a trade-off between all these methods.
It is only slightly slower than~$\methodatck$, and its convergence continues to frequently improve the incumbent solution.
The returned solution is feasible and has a high objective value.

\begin{figure}[!ht]
    \centering
    \includegraphics[width=0.9\linewidth]{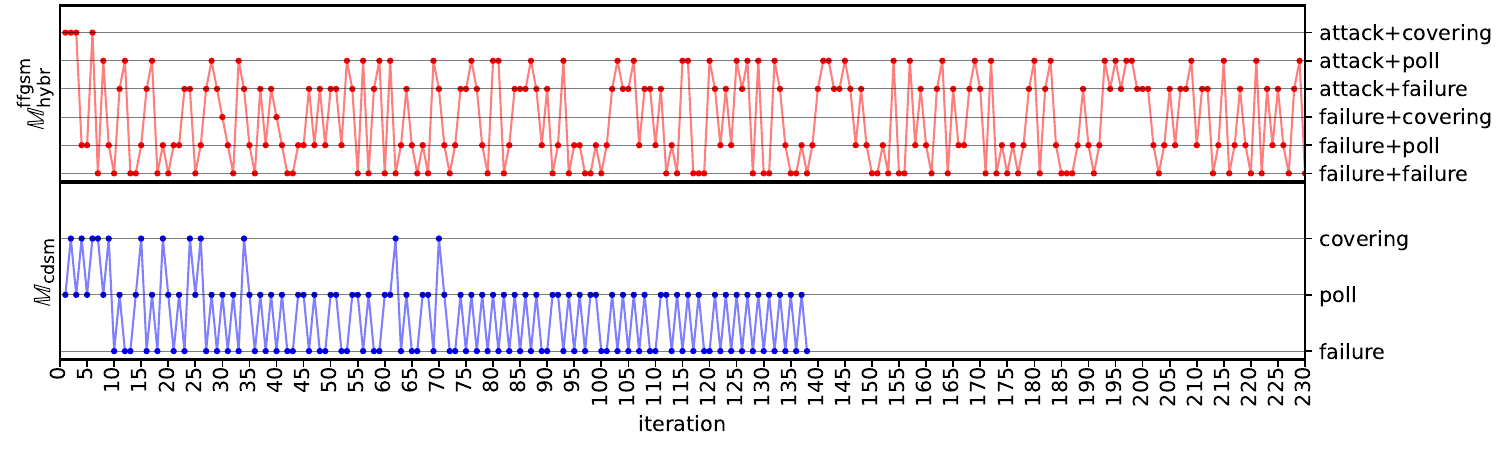}
    \caption{\nameref{experiment:contribution_steps} in the Counterfactual Warcraft Maps problem from Section~\ref{section:numerical/warcraft}.}
    \label{figure:warcraft/iterations_results}
\end{figure}

\nameref{experiment:contribution_steps}, illustrated in Figure~\ref{figure:warcraft/iterations_results}, shows that the \pollstep step in the \cdsm-based methods struggles.
The constraint that~$x$ must remain close to the unit sphere is challenging for the \pollstep step, since a large \pollstep radius has little chance to preserve the feasibility of this constraint while a small radius provides little improvements.
This explains the failure-success pattern that we observe for the~$\methodcdsm$ method.
The \pollstep in the~$\methodhybr$ method follows a similar pattern, but this is mitigated by the~$\attackstep$ step that succeeds roughly half of the time.
This, combined with the earlier observation that the~$\methodatck$ method alone fails on a worse solution, suggests that the \cdsm components of the~$\methodhybr$ method reposition the current incumbent in more favourable areas of the input space when the \attackstep fails, so that new successful attacks identify new ascent directions.

\begin{figure}[!ht]
    \centering
    \includegraphics[width=0.9\linewidth]{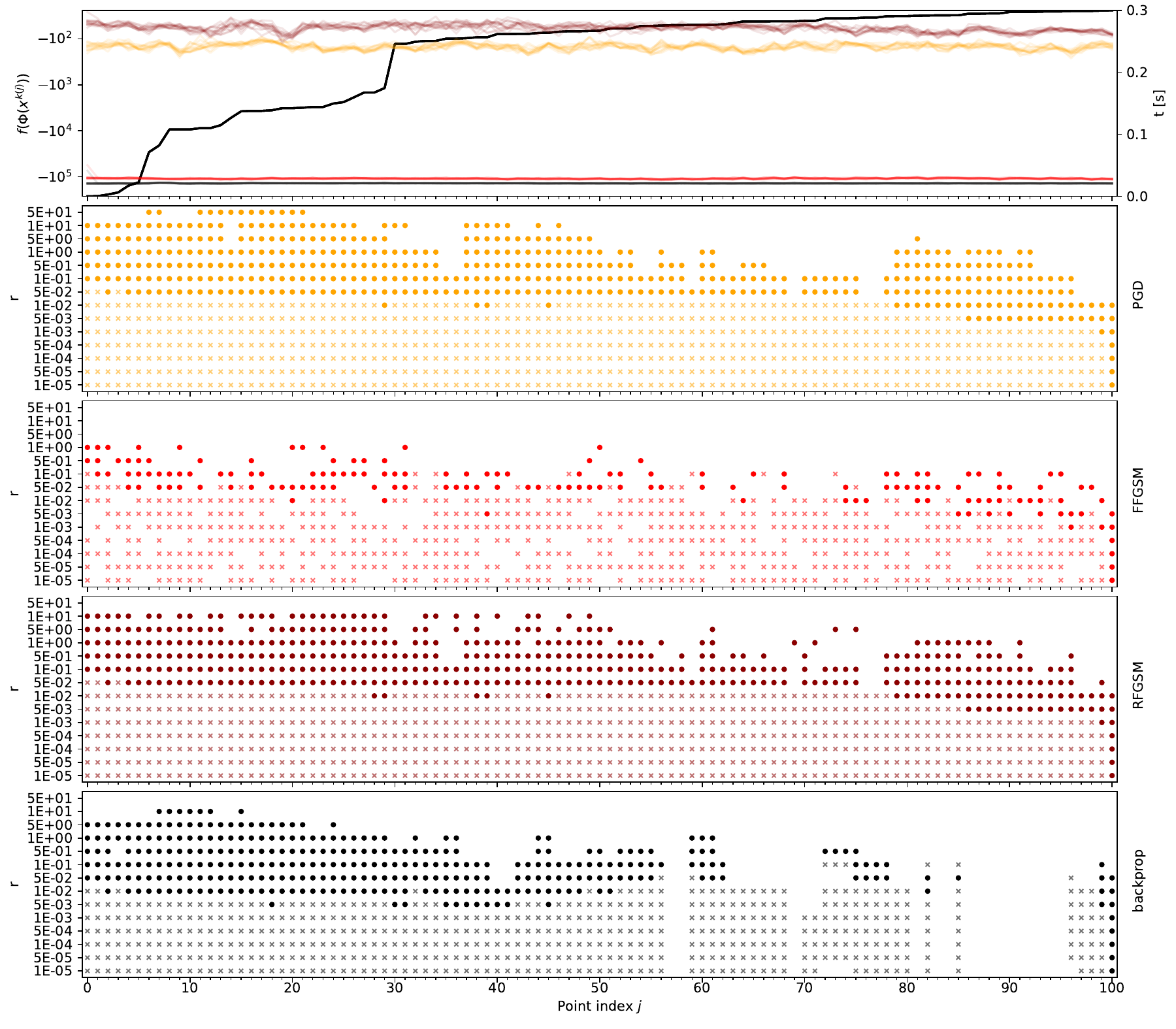}
    \caption{\nameref{experiment:potential_attack} in the Counterfactual Warcraft Maps problem from Section~\ref{section:numerical/warcraft}.}
    \label{figure:warcraft/attacks_analysis}
\end{figure}

\nameref{experiment:potential_attack} corroborates these observations.
Figure~\ref{figure:warcraft/attacks_analysis} shows that the \attackstep operator reliably identifies ascent directions all along the optimization process, while the \backprop-based operator sometimes fails (especially at~$j \approx 90$).
The \rfgsm and \pgd variants achieve comparable success rates across the sample, yet~\rfgsm is substantially faster.

\subsection{Application: Bio-Diesel Production}
\label{section:numerical/bio_pinn}

\newcommand{\specieTG}{\mathrm{TG}}
\newcommand{\specieDG}{\mathrm{DG}}
\newcommand{\specieMG}{\mathrm{MG}}
\newcommand{\specieG}{\mathrm{G}}
\newcommand{\specieME}{\mathrm{ME}}

This problem is constructed from the chemical engineering study of bio-diesel production from~\cite{BiBoBl24PINNbiodiesel}.
The process they study involves five chemical species: the desired methyl ester ($\specieME$), and four intermediate reactants: triglycerides ($\specieTG$), diglycerides ($\specieDG$), monoglycerides ($\specieMG$), and glycerol ($\specieG$).
For each species, we denote by~$[\cdot](t,Q)$ its concentration (in moles per litre) in the reactor tank after a reaction of duration~$t$ (in seconds) under constant incoming heat flux with power~$Q$ (in Watts).
We also denote by~$T(t,Q)$ the reactor temperature (in Celsius degrees) under these conditions.

We consider~$K \defequal 40$ reactor tanks, all operating under the same initial conditions.
We optimize over~$(t_k,Q_k)_{k=1}^{K} \in (\Rbb_+^2)^K$ to maximize the production of~$\specieME$, measured in each tank via the final purity ratio of the mixture it contains.
To prevent evaporation of~$\specieME$ (with boiling point at~$65^\circ\Crm$), we require for~$T(s,Q_k) \leq 65$ for all~$s \in [0,t_k]$ and all~$k \in \llb1,K\rrb$.
Moreover, we enforce a total energy budget of~$300K$ currency units.
For all~$k \in \llb1,K\rrb$, the energetic cost of the tank~$k$ is~$w_kQ_kt_k$ currency units, where~$w_k \defequal 1+\frac{3}{4}\sin(2\pi\frac{k-1}{K})$ (in currency per Joule).
The resulting problem reads as
\begin{equation*}
    \problemoptim
        {\maximize}
        {(t_k)_{k=1}^{K} \in \Rbb_+^K \\ (Q_k)_{k=1}^{K} \in \Rbb_+^K}
        {\dfrac{1}{K} \dsum{k=1}{K} \dfrac{[\specieME(t_k,Q_k)]}{[\specieTG](t_k,Q_k) + [\specieDG](t_k,Q_k) + [\specieMG](t_k,Q_k) + [\specieG](t_k,Q_k) + [\specieME](t_k,Q_k)}}
        {T(s,Q_k) \leq 65,~ \forall s \in [0,t_k], ~\forall k \in \llb1,K\rrb, \quad \text{and} \quad \dsum{k=1}{K} w_kQ_kt_k \leq 300K.}
\end{equation*}
Note that the objective function is symmetric with respect to the~$K$ reactor tanks.
Only the energy constraint breaks the symmetry on the problem thanks to the nonconstant weights.

We solve this problem by a simulation-based approach where we simulate the chemical process with the \pinn trained in~\cite{BiBoBl24PINNbiodiesel}, available at~\url{https://github.com/chaos-polymtl/bio-pinn}.
The \pinn predicts, for any given reaction time~$t$ and heating power~$Q$, the concentrations of the species involved in the reaction and the temperature at the end of the process.
The model reads as a function from the input space~$\Ibb \defequal \Rbb_+^2$ to the predictions space~$\Pbb \defequal \Rbb^6$ given by
\begin{equation*}
    \fct
        {\pinn}
        {(t,Q)}
        {\Ibb}
        {\begin{bmatrix} [\specieTG](t,Q) ~~ [\specieDG](t,Q) ~~ [\specieMG](t,Q) ~~ [\specieG](t,Q) ~~ [\specieME](t,Q) ~~ T(t,Q) \end{bmatrix}.}
        {\Pbb}
\end{equation*}
Since the \pinn only provides predictions at the queried time~$t$, we estimate intermediate values by discretizing~$[0,t]$ into~$N+1$ steps (with~$N \defequal 10$) and evaluating the~$\pinn$ at each step.
We thus define
\begin{equation*}
    \fct
        {\Psi}
        {(t,Q)}
        {\Ibb}
        {\left( \pinn\!\left(\tfrac{it}{N}, Q\right) \right)_{i=0}^{N}.}
        {\Pbb^{N+1}}
\end{equation*}
To ensure meaningful predictions, we impose additional constraints.
First, the input domain is bounded by~$0 \leq t \leq 600$ and~$0 \leq Q \leq 12$, following~\cite{BiBoBl24PINNbiodiesel}, to remain near the training regime of the model.
Second, all concentrations are required to be nonnegative.
While this holds physically, the constraint prevents occasional violations due to modelling imperfections by the \pinn.

The problem is expressed as an instance of Problem~\eqref{problem:P} by defining the variables space~$\Xbb \defequal \Ibb^K$, the outputs space~$\Obb \defequal (\Pbb^{N+1})^K$, and the neural mapping
\begin{equation*}
    \fct
        {\Phi}
        {x \defequal (t_k, Q_k)_{k=1}^{K}}
        {\Xbb}
        {\left( \left(\Psi(t_k, Q_k)\right)_{k=1}^{K}, x \right).}
        {\Obb \times \Xbb}
\end{equation*}
Hereafter, for all~$z \in \Obb$, all~$k \in \llb1,K\rrb$, all~$i \in \llb0,N\rrb$ and all~$u \in \llb1,6\rrb$, we denote by~$z_{(k,i,u)}$ the~$u$-th prediction at the~$i$-th time step in the batch~$k$.
Then, we define the objective function~$f$ as
\begin{equation*}
    \fct
        {f}
        {(z,x)}
        {\Obb \times \Xbb}
        {\dfrac{1}{K} \dsum{k=1}{K} \dfrac{ z_{(k,N,5)} }{ z_{(k,N,1)} + z_{(k,N,2)} + z_{(k,N,3)} + z_{(k,N,4)} + z_{(k,N,5)} },}
        {\Rbb}
\end{equation*}
and our constraint function by
\begin{equation*}
    \fct
        {c}
        {\left(z, x \defequal ((t_k, Q_k)_{k=1}^{K})\right)}
        {\Obb \times \Xbb}
        {\left(
            \dsum{k=1}{K} w_kQ_kt_k - 300K,~
            \begin{bmatrix} -t_k \\ t_k-600 \\ -Q_k \\ Q_k-12 \end{bmatrix}_{k=1}^{K},~
            \begin{bmatrix} -z_{k,i,1} \\ -z_{k,i,2} \\ -z_{k,i,3} \\ -z_{k,i,4} \\ -z_{k,i,5} \\ z_{k,i,6}-65 \end{bmatrix}_{k=0,~ i=0}^{K,N}
            \right).}
        {\Rbb \times (\Rbb^{4})^K \times \Rbb^{6K(N+1)}}
\end{equation*}
Then, for all~$x \defequal (t_k, Q_k)_{k=1}^{K} \in \Xbb$,
\begin{equation*}
    f(\Phi(x)) = \dfrac{1}{K} \dsum{k=1}{K} \dfrac{ [\specieME]\left(t_k,Q_k\right) }{ [\specieTG]\left(t_k,Q_k\right) + [\specieDG]\left(t_k,Q_k\right) + [\specieMG]\left(t_k,Q_k\right) + [\specieG]\left(t_k,Q_k\right) + [\specieME]\left(t_k,Q_k\right) }
\end{equation*}
recovers the objective of the problem, and the constraints are enforced at the discretized times since
\begin{equation*}
    c(\Phi(x)) \leq 0 \iff
    \begin{cases}
        \sum_{k=1}^{K} w_kQ_kt_k \leq 300K,
        \\
        0 \leq t_k \leq 600, & \forall k \in \llb1,K\rrb,
        \\
        0 \leq Q_k \leq 12, & \forall k \in \llb1,K\rrb,
        \\
        [\Srm](\frac{it_k}{N}, Q_k) \geq 0, &
            \forall k \in \llb1,K\rrb,~
            \forall i \in \llb0,N\rrb,~
            \forall \Srm \in \{\specieME, \specieTG, \specieDG, \specieMG, \specieG\},
        \\
        T(\frac{it_k}{N}, Q_k) \leq 65, &
            \forall k \in \llb1,K\rrb,~
            \forall i \in \llb0,N\rrb,
    \end{cases}
\end{equation*}
We consider~$x^0 \defequal (60, 2)_{k=1}^{K}$ and~$\lambda \defequal 10$.
Note that~$f$ has active subspaces, since for all~$(z,x) \in \Obb \times \Xbb$, the value~$f(z,x)$ is independent to~$x$ and all~$z_{(k,i,u)}$ with either~$i < N$ or~$(i,u) = (N,6)$, so we adapt the \attackstep operator as in Remark~\ref{remark:active_subspaces}.
All methods return solutions that slightly break a few constraints enforcing positive concentrations at the early time steps, but the maximal violation ($0.005$ for the~$\methodcdsm$ method and~$0.012$ for all other methods) is below the experimental uncertainty in the training data of the \pinn.
The energy budget is only exceeded in two runs of the~$\methodatck$ method, in both case by~$0.002$.

\begin{figure}[!ht]
    \centering
    \includegraphics[width=0.9\linewidth]{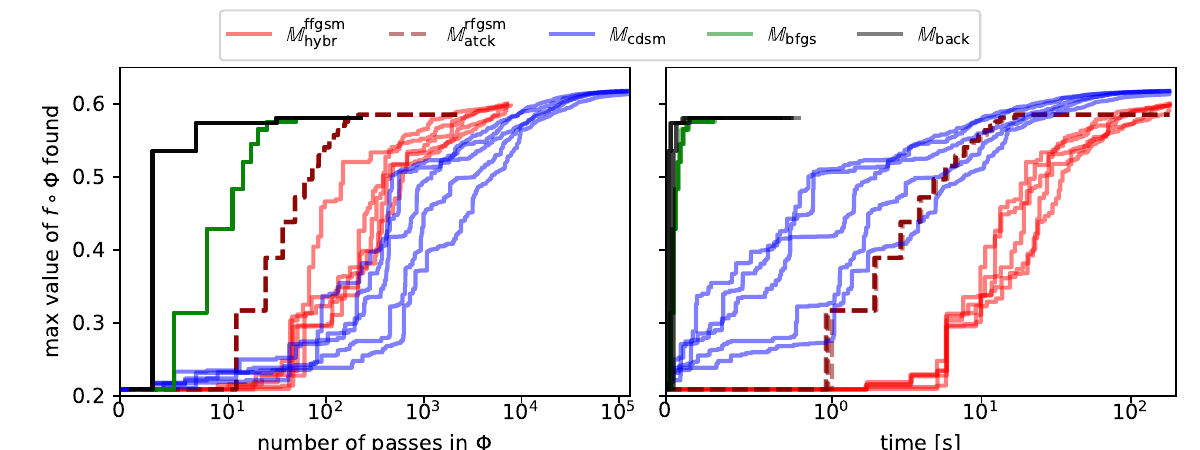}
    \caption{\nameref{experiment:performance_algo} in the Bio-Diesel Production problem from Section~\ref{section:numerical/bio_pinn}.}
    \label{figure:bio_pinn/results_nb_pts}
\end{figure}

\nameref{experiment:performance_algo}, shown in Figure~\ref{figure:bio_pinn/results_nb_pts}, highlights that all gradient-based methods converge to the same objective value.
The~$\methodbprp$ and~$\methodbfgs$ methods actually both return the point~$x_\sharp \approx (60.14, 4.99)_{k=1}^{K}$ (remark that the~$K$ couples in~$x_\sharp$ are all equal).
This results from the symmetry in the objective function.
The~$\methodatck$ method returns a point close to~$x_\sharp$ despite the randomness in the \rfgsm algorithm.
In contrast, the~$\methodcdsm$ and~$\methodhybr$ methods converge to points with no symmetries, and that significantly differ from~$x_\sharp$.
Finally, Figure~\ref{figure:bio_pinn/results_nb_pts} shows that the~$\methodatck$ and~$\methodhybr$ methods are slower than the~$\methodcdsm$ method despite calling the \nn~$\Phi$ less often.
This has little theoretical interpretation, so we suspect that the implementations of the \attackstep algorithms struggle with the separable structure of~$\Phi$.

\begin{figure}[!ht]
    \centering
    \includegraphics[width=0.9\linewidth]{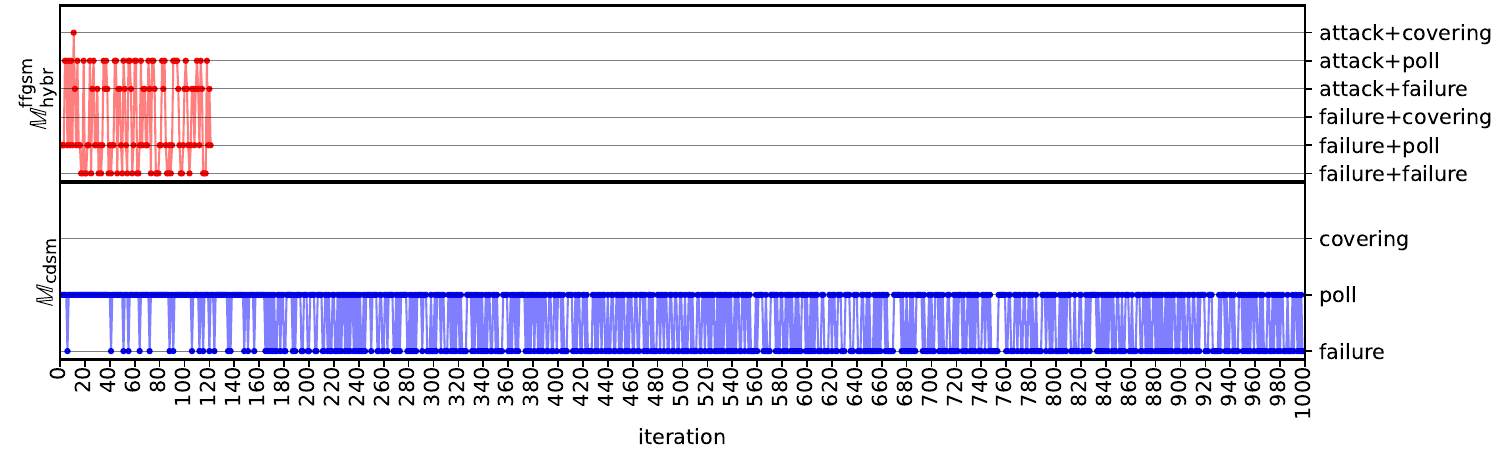}
    \caption{\nameref{experiment:contribution_steps} in the Bio-Diesel Production problem from Section~\ref{section:numerical/bio_pinn}.}
    \label{figure:bio_pinn/iterations_results}
\end{figure}

Now we focus to \nameref{experiment:contribution_steps}, illustrated in Figure~\ref{figure:bio_pinn/iterations_results}.
The~$\methodhybr$ method stops (as it exceeds the allowed runtime) after~$10$ times fewer iterations than the~$\methodcdsm$ method.
As already discussed, this results from the slow attack algorithm.
The~$\methodcdsm$ method faces numerous \pollstep step failures, presumably because the different scales of the variables make the feasible set difficult to explore.
In the~$\methodhybr$ method, the \pollstep step remains affected by this behaviour, yet it consistently repositions the incumbent solution so that the \attackstep step often succeeds.

\begin{figure}[!ht]
    \centering
    \includegraphics[width=0.9\linewidth]{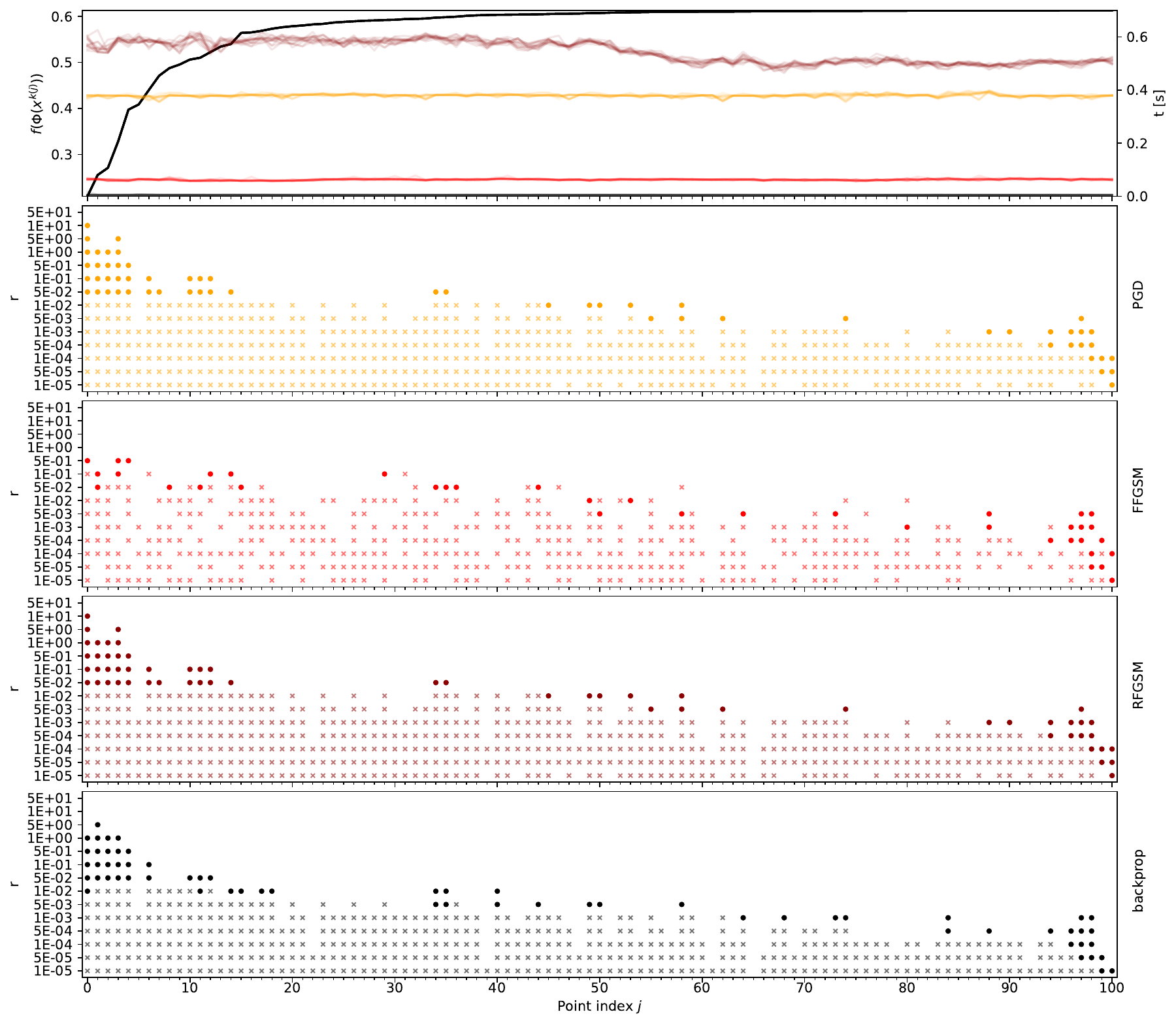}
    \caption{\nameref{experiment:potential_attack} in the Bio-Diesel Production problem from Section~\ref{section:numerical/bio_pinn}.}
    \label{figure:bio_pinn/attacks_analysis}
\end{figure}

\nameref{experiment:potential_attack}, shown in Figure~\ref{figure:bio_pinn/attacks_analysis}, further confirms that the \attackstep operator is potent in early optimization or at low radius in late optimization.
The performances of the \pgd and \rfgsm algorithm equal those of the backpropagation, while \ffgsm appears less robust.

\section{Conclusions and Future Work}
\label{section:conclusion}

To conclude this paper, we discuss some strengths and limitations of our hybrid method in Section~\ref{section:conclusion/discussion}, and we outline promising avenues for future research in Section~\ref{section:conclusion/perspectives}.

\subsection{Strengths and Limitations of our Hybrid Method}
\label{section:conclusion/discussion}

The theoretical analysis of Algorithm~\ref{algorithm:hybrid}, given by Theorem~\ref{theorem:convergence_hybrid}, is limited to an asymptotic convergence.
It is difficult to conduct a non-asymptotic analysis under the broad Assumption~\ref{assumption:problem}.
Such analyses are scarce in the \dfo literature, unless significant additional assumptions on the objective and constraints functions are added.
Nevertheless, a strength of Theorem~\ref{theorem:convergence_hybrid} is that Algorithm~\ref{algorithm:hybrid} is guaranteed to identify a local solution to Problem~\eqref{problem:P}, even in hard instances.
Another positive trait of the \coveringstep step is that Algorithm~\ref{algorithm:hybrid} scans a dense set of points in a neighbourhood around that solution, which eventually gives a precise understanding of the landscape of Problem~\eqref{problem:P} near that solution.

Our hybrid method enjoys a twofold flexibility.
The core components of Algorithm~\ref{algorithm:hybrid} (the \attackstep operator and the \cdsm routine) are largely independent.
Each may be chosen depending on the problem at hands.
This flexibility allows our hybrid method to be adapted to a wide range of contexts, but a downside is a concern on the choice of the most suitable components.
Let us discuss some guideline about how to define each component of our hybrid method depending on the context.

First, the \attackstep operator may be defined from several algorithms for \nn attacks.
An ideal choice depends on whether~$\Phi$ is given as a glass-box \nn.
Indeed, if~$\Phi$ is a glass-box \nn, we recommend using an algorithm for \nn attacks that exploits this explicit structure, such as those in the \abcrown solver.
If~$\Phi$ is not a glass-box \nn, then we suggest using algorithms based on \fgsm (such as \fgsm itself, or \rfgsm or \ffgsm) since they are fast and our numerical experiments hint that more accurate algorithm such as \pgd yield negligible difference in the performance of the \attackstep operator.
However, it is also worth trying to use backpropagations directly to compute ascent directions.
Our experiments show that it is sometimes prone to early failure, but nonetheless it is very fast and works well on some of our numerical experiments.

Second, in Algorithm~\ref{algorithm:hybrid}, we choose to hybridize a generic \attackstep operator with the \cdsm from \dfo.
Others algorithms could be substituted for \cdsm, but unfortunately, to the best of our knowledge, no consensus currently exists in the literature about \dfo regarding the most appropriate choice for a given problem.
The \cdsm is suited for cases where either~$\Phi$ is not a glass-box \nn or~$f$ or~$c$ are generic nonlinear functions, so the structure of the problem is limited.
However, in other contexts, we suggest considering hybridizing the \attackstep step with other algorithms instead.
A \dfo method is likely not the most efficient approach when~$\Phi$ is a glass-box \nn and either~$f$ and~$c$ are simple enough.
For example, when methods from the literature (see Section~\ref{section:related_literature}) about optimization through trained \nn are applicable (that is, when~$\Phi$ is a \relu-based glass-box \nn and~$f$ is linear and~$c$ yields a polyhedral feasible set~$F$), they presumably should be preferred.
Similarly, when~$f$ is nonlinear but simple enough (for example, quadratic) and~$c$ yields a polyhedral feasible set, then it is presumably more efficient to adapt these methods than to consider a \dfo method.
When usable, methods such as \bfgs may also be considered (possibly by replacing the composite gradient by an \attackstep direction).

\subsection{Perspectives for Future Research}
\label{section:conclusion/perspectives}

Numerous avenues for research stem from our work, either on the theoretical and practical sides.

First, we may strengthen the \attackstep component itself, both in the choice of ascent directions for~$f$ and in the way attacks are computed.
Beyond gradients, Remark~\ref{remark:which_ascent_direction_to_attack} suggests using alternative ascent proxies such as simplex or finite-difference gradients when~$f$ lacks smoothness.
Adapting our analysis to these settings is straightforward, and the literature on simplex gradients and simplex Hessians~\cite{HaJBCh20Hessian,HaJBPl23PositiveBases,HaJBPl20ErrorsBoundsSimplexGrad} offers principled constructions that could yield more reliable attack directions.
On the numerical side, our experiments indicate that speed often matters more than ultimate attack accuracy.
This makes \fgsm an appealing default option, although this choice is likely problem-dependent.
In particular, \fgsm is likely not sufficient in cases where the structure of~$\Phi$ makes it difficult to compute successful attacks.
More broadly, we may also develop attack solvers that natively handle input constraints and possess guarantees of success whenever an attack exists.
Such constrained attack would close the gap with Assumption~\ref{assumption:attack_guarantees_ascent_directions} and remove the heuristic workaround from Section~\ref{section:optimization_attack/use_in_optim/heuristic_tool}, and would also be broadly useful beyond our context.
Alternatively, we may also reconsider our workaround penalizing constraints in the directional attack problem, to ensure that it remains exact.
Some techniques such as exact penalization theory or augmented Lagrangian method may relate to such research.

Second, we plan to broaden the scope of problems addressed, with a particular emphasis on regularity and constraints.
The \cdsm copes with possible discontinuities~\cite{AuBoBo24Covering}, and the assumptions underlying its convergence are \textit{tight}.
This makes the extension of Algorithm~\ref{algorithm:hybrid} to discontinuous~$f$,~$c$, or~$\Phi$ natural, provided the \attackstep follows the ascent-proxy adaptations discussed above.
On the constraints side, we may replace strict feasibility or global relaxation with mature mechanisms from \dfo such as the progressive barrier~\cite{AuDe09a}.
This could improve efficiency by allowing controlled, temporary infeasibility, which is common in successful \dfo practice.

Third, we expect benefits from enhancing the \dfo engine that surrounds the \attackstep step in our hybrid method.
Within \cdsm, precision and speed can be improved by established tools, such as local quadratic models~\cite{CoLed2011,VeKa2017} to inform \pollstep directions, Bayesian strategies~\cite{ZhZi21GlobalBayesian} to steer the \searchstep step, and variable neighbourhood rules~\cite{HaMl01a,MlDrKoVuCa2008}.
These improvements are independent of the attack component and can be integrated without altering convergence guarantees, offering routes to better late-stage refinement.
We may also consider a hybrid algorithm that drives the \attackstep radius at each iteration according to other strategies, such as trust region or momentum-based approaches.

Fourth, we see opportunities to leverage problem structure more explicitly.
When~$\Phi$ is a glass-box \nn, hybridizing the \attackstep with methods tailored to optimization through trained \nns (see Section~\ref{section:related_literature}) should be preferred to structure-agnostic \dfo methods.
Moreover, the composite form invites alternative formulations, e.g., introducing an auxiliary variable~$y$ representing the output of~$\Phi(x)$, with the constraint that~$y = \Phi(x)$, and optimizing~$f(y)$ under~$c(y) \leq 0$.
These viewpoints connect naturally with partitioned~\cite{AuBoBo24Partitioned} and parametric optimization~\cite{St18Parametric}, and with optimization on manifolds~\cite{boumal23Manifold}.
Combining attack-based steps with such dedicated methods may unlock further gains of performance.
For example, we may allow to temporarily break the constraint~$y = \Phi(x)$, so that whenever a relevant value for~$y$ is found, a targeted attack of~$\Phi$ with target~$y$ may be used to modify~$x$ accordingly.

Finally, a broader experimental campaign across additional architectures (e.g., physics-informed models~\cite{CuDCGiRoRaPi22PINNsurvey,ShFeChZhCaJi25PINNSurvey} and digital twins~\cite{TaXiQiChJi22DigitalTwin}) would clarify when lightweight attacks suffice and when stronger attacks are warranted.
Although the above fields are active, we identified few pre-trained and publicly available \nn surrogates that suit our needs.
Similarly, current solvers for \nn attacks are not compatible with \nn that do not provide access to backpropagation.
Then, our current experiments do not show how well an \attackstep operator relying on such attack algorithms performs in practice.
If future solvers for \nn attacks implement algorithms working in black-box \nn access, an experimental campaign should test them in this regime.

\printbibliography

\end{document}